\input epsf.sty

\def\cD{{\cal D}}
\def\cF{{\cal F}}
\def\cG{{\cal G}}
\def\cI{{\cal I}}
\def\cV{{\cal V}}
\def\cX{{\cal X}}
\def\cY{{\cal Y}}

\def\bC{{\bf C}}
\def\bC{{\bf C}}
\def\bC{{\bf C}}
\def\bN{{\bf N}}
\def\bP{{\bf P}}
\def\bR{{\bf R}}
\def\bT{{\bf T}}
\def\bZ{{\bf Z}}

\def\sqr#1#2{{\vcenter{\hrule height.#2pt              
     \hbox{\vrule width.#2pt height#1pt\kern#1pt
     \vrule width.#2pt}
     \hrule height.#2pt}}}
\def\square{\mathchoice\sqr{5.5}4\sqr{5.0}4\sqr{4.8}3\sqr{4.8}3}
\def\qed{\hskip4pt plus1fill\ $\square$\par\medbreak}

\magnification\magstep1
\centerline{\bf Dynamics of Rational Surface Automorphisms:}

\centerline{\bf Linear Fractional Recurrences}

\bigskip
\centerline{ Eric Bedford\footnote*{Research supported in part by the NSF.}  and Kyounghee Kim}
\bigskip
\bigskip\noindent{\bf \S0.   Introduction.  }  Here we discuss automorphisms (biholomorphic maps) of compact, projective surfaces with positive entropy.  Cantat [C1] has shown that the only (minimal dynamical system model) possibilities occur for tori and $K3$  (and certain of their quotients), and rational surfaces.  $K3$ surfaces have been studied by Cantat [C2] and McMullen [M1].  Here we consider the family of birational maps of the plane which are defined by 
$$f_{a,b}: (x,y)\mapsto \left(y,{y+a\over x+b}\right),\eqno(0.1)$$
and which provide an interesting source of automorphisms of rational surfaces.  
The maps $f_{a,b}$ form part of the family of so-called linear fractional recurrences, which were studied from the point of view of degree growth and periodicity in [BK].  We let  $\cV=\{(a,b)\in\bC^2\}$ be the space of parameters for this family, and we define 
$$\eqalign{q&=(-a,0), \ p=(-b,-a),\cr
\cV_n&=\{(a,b)\in\cV:f_{a,b}^jq\ne p {\rm \ for\ }0\le j<n,{\rm \ and\ }f_{a,b}^nq=p\}.}\eqno(0.2)$$  
In [BK] we showed that $f_{a,b}$ is birationally conjugate to an automorphism  of a compact, complex surface $\cX_{a,b}$ if and only if $(a,b)\in\cV_n$ for some $n\ge0$.  The surface $\cX_{a,b}$ is obtained by blowing up the projective plane $\bP^2$ at the $n+3$ points $e_1=[0:1:0]$, $e_2=[0:0:1]$, and $f^jq$, $0\le j\le n$.  The dimension of $Pic(\cX_{a,b})$ is thus $n+4$, and the action of $f_{a,b}^*$ on $Pic(\cX_{a,b})$ is the same for all $(a,b)\in\cV_n$: its characteristic polynomial is
$$\chi_n(x):=-1 + x^2+x^3-x^{1+n} - x^{2+n} + x^{4+n}.\eqno(0.3)$$
This polynomial arises in the growth of Coxeter groups (see [F, p.\ 483]).  When $n\ge7$, $\chi_n$ has a root $\lambda_n>1$, which is the unique root with modulus greater than one; and the entropy of $f_{a,b}$ is $\log\lambda_n>0$.   If $\psi_n(x)$ denotes the minimal polynomial of $\lambda_n$, then we may factor $\chi_n(x) = C_n(x) \psi_n(x)$, where $C_n$ is a product of cyclotomic polynomials.   The factorization of $C_n$ into cyclotomic polynomials is given explicitly in Theorems 3.3 and 3.5.

McMullen [M2] noted that the maps $f_{a,b}$ provide representations of the Coxeter elements of certain Weyl groups, and he gave a construction of automorphisms which represent these Coxeter elements.  By Theorem 8.7,  the maps $f_{a,b}$ give essentially all possible rational surface automorphisms which represent these Coxeter elements.  

The approach of this paper is to focus on the maps $f_{a,b}$ that have invariant curves.  We obtain the following in Theorem 4.2:
\proclaim Theorem A.  There are rational surface automorphisms $f_{a,b}$ which have positive entropy but have no invariant curves.  

\noindent We show in \S4 that this gives a counterexample of a conjecture/question of Gizatullin, Harbourne, and McMullen.

We are also interested in the dynamical behavior of the maps with positive entropy.  We say that a connected,  invariant open set ${\cal D}$ is a rotation domain if the restriction $f|_{{\cal D}}$ is conjugate to a rotation of infinite order.  (For mappings in dimension one, this corresponds to having a Siegel disk.)  For a rotation domain ${\cal D}$,  the closure of a generic orbit will be a smooth real torus ${\bf T}^d$ with $d$ equal to 1 or 2, and we say that $d$ is the rank of the rotation domain.

\proclaim Theorem B.    Let $f_{a,b}$ be a mapping of the form (0.1) which is equivalent to an automorphism with positive entropy.  If there is an $f_{a,b}$-invariant curve, then exactly one of the following occurs:
\item{(i)}  $f$ has a rank 1 rotation domain centered at one of the two fixed points.
\item{(ii)} $a,b\in\bR$, and the restriction $f_R$ of $f$ to the real points ${\cal X}_R$ has the same entropy as $f$.  Further the (unique) invariant measure of maximal entropy is supported on a subset of $\cX_R$ of zero area.

As $n$ increases,  the maps with invariant cubics appear to  form a very small fraction of the total number of maps in $\cV_n$.  It will be interesting to explore the dynamical properties of these other maps.

This paper is organized as follows.  In \S1 we show that if $S_{a,b}$ is a curve which is invariant under a map $f_{a,b}$ as in (0.1), then $S_{a,b}$ is cubic.  We give the precise form of the cubic curve $S_{a,b}$ (if it exists) in \S2.  In \S2, we also show that these parameter values $(a,b)$ corresponding to invariant curves, correspond to $(a,b)\in\Gamma$, where $\Gamma$ is the union of three specific curves in parameter space.  In \S3 we show that there is a close connection between the zeros of $\chi_n$ and the parameter values for which $f_{a,b}$ is an automorphism.   In particular, we find the cyclotomic part $C_n$ of $\chi_n$.  In \S4 we show the existence of counterexamples to the Gizatullin/Harbourne/McMullen conjecture.  And we prove Theorem 4.2, which also contains Theorem A above.  We discuss the existence of rotation domains in \S5.  In \S6 we show that the real maps with invariant curves have entropy equal to $\log\lambda_n$.  In \S7 we prove Theorem B.  In \S8 we discuss the connection with representations of the Coxeter element.

Acknowledgement.  We thank Serge Cantat and Jeff Diller for explaining some of this material to us and giving helpful suggestions on this paper.  We also thank the anonymous referees for their helpful comments.

\bigskip\noindent{\bf \S1.  Invariant Curves.  }     We will write our maps in projective coordinates $x=[x_0:x_1:x_2]=[1:x:y]$.  We define $\alpha=\gamma=(a,0,1)$ and $\beta=(b,1,0)$, so $f$ is written
$$f_{a,b}:\ [x_0:x_1:x_2]\mapsto [x_0\beta\cdot x:x_2\beta\cdot x:x_0\alpha\cdot x].$$ 
The exceptional curves for the map $f$ are given by the lines $\Sigma_0=\{x_0=0\}$, $\Sigma_\beta=\{\beta\cdot x=0\}$, and $\Sigma_\gamma=\{\gamma\cdot x=0\}$.   The indeterminacy locus $\cI(f)=\{e_2,e_1,p\}$ consists of the vertices of the triangle $\Sigma_0\Sigma_\gamma\Sigma_\beta$.    Let $\pi:\cY\to\bP^2$ be the complex manifold obtained by blowing up $e_1$ and $e_2$, and let the exceptional fibers be denoted $E_1$ and $E_2$.  By $\Sigma_0$, $\Sigma_\beta$ and $\Sigma_\gamma$ we denote the strict transforms in $\cY$.  Let $f_\cY:\cY\to\cY$ be the induced birational map.  Then the exceptional locus is $\Sigma_\gamma$, and the indeterminacy loci are $\cI(f_\cY)=\{p\}$ and $\cI(f_\cY^{-1})=\{q\}$.  In particular, $f_\cY:\Sigma_\beta\to E_2\to\Sigma_0\to E_1\to \Sigma_B=\{x_2=0\}.$  By curve, we mean an algebraic set of pure dimension 1, which may or may not be irreducible or connected.
We say that an algebraic curve $S$ is invariant if the closure of $f(S-\cI)$ is equal to $S$.   We define the cubic polynomial $j_f:=x_0(\beta\cdot x)(\gamma\cdot x)$, so $\{j_f=0\}$ is the exceptional locus for $f$.  For a homogeneous polynomial  $h$ we consider the condition that there exists $t\in\bC^*$ such that
$$h\circ f=t\cdot j_f\cdot h.\eqno(1.1)$$

\proclaim Proposition 1.1.  Suppose that $(a,b)\notin\bigcup\cV_n$, and $S$ is an $f_{a,b}$-invariant curve.  Then  $S$ is a cubic containing $e_1,e_2$, as well as $f^jq,f^{-j}p$ for all $j\ge0$.  Further, (1.1) holds for $S$.

\noindent{\it Proof. }   Let us pass to $f_\cY$, and let $S$ denote its strict transform inside $\cY$.  Since $(a,b)\notin\bigcup\cV_n$, the backward orbit $\{f^{-n}p:n\ge1\}$ is an infinite set which is disjoint from the indeterminacy locus of $\cI(f_\cY^{-1})$.  It follows that $S$ cannot be singular at $p$  (cf.\ Lemma 3.2 of [DJS]).  Let $\mu$ denote the degree of $S$, and let $\mu_1=S\cdot E_1$, $\mu_0=S\cdot\Sigma_0$, and $\mu_2=S\cdot E_2$.  It follows that $\mu=\mu_2+\mu_0+\mu_1$.  Further, since $f_\cY:E_2\to\Sigma_0\to E_1$, we must have $\mu_2=\mu_0=\mu_1$.  Thus $\mu$ must be divisible by 3, and $\mu/3=\mu_2=\mu_0=\mu_1$.  Now, since $S$ is nonsingular at $p=\Sigma_\beta\cap\Sigma_\gamma$, it must be transversal to either $\Sigma_\gamma$ or $\Sigma_\beta$.  

Let us suppose first that $S$ is transverse to $\Sigma_\gamma$ at $p$.  Then we have $\mu=\mu_2+S\cdot\Sigma_\gamma$.  Thus $S$ intersects $\Sigma_\gamma-\{p\}$ with multiplicity $\mu-(\mu/3)-1$.  If $\mu>3$, then this number is at least 2.  Now $\Sigma_\gamma$ is exceptional,  $f_\cY$ is regular on $\Sigma_\gamma-\{p\}$, and $f_\cY(\Sigma_\gamma-\{p\})=q$.  We conclude that $S$ is singular at $q$.  This is not possible by  Lemma 3.2 of [DJS] since $q$ is indeterminate for $f^{-1}$.  This is may also be seen because since $(a,b)\notin\bigcup\cV_n$, it follows that $f^nq$ is an infinite orbit disjoint from the indeterminacy point $p$, which is a contradiction since $S$ can have only finitely many singular points.

Finally, suppose that $S$ is transversal to $\Sigma_\beta$.  We have $\mu=S\cdot\Sigma_\beta+\mu_2$, and by transversality, this means that $S$ intersects $\Sigma_\beta-\{p\}$ with multiplicity ${2\mu\over 3}-1$.  On the other hand, $f_\cY$ is regular on $\Sigma_\beta-\{p\}$, and $\Sigma_\beta-\{p\}\to E_2$.  Thus the multiplicity of intersection of $S$ with $\Sigma_\beta-\{p\}$ must equal the multiplicity of intersection with $E_2$, but this is not consistent with the formulas unless $\mu=3$. \qed

The following was motivated by [DJS]:  
\proclaim Theorem 1.2.  Suppose that $(a,b)\in\cV_n$ for some $n\ge11$.  If $S$ is an invariant algebraic curve, then the degree of $S$ is 3, and (1.1) holds.

\noindent{\it Proof. }  Let $\cX$ be the manifold $\pi:\cX\to\bP^2$ obtained by blowing up $e_1,e_2,q,fq,\dots,f^nq=p$, and denote the blowup fibers by $E_1,E_2,Q,fQ,\dots,f^nQ=P$.   Suppose that $S$ is an invariant curve of degree $m$.   By $S$, $\Sigma_0$, etc., we denote the strict transforms of these curves inside $\cX$.  Let $f_\cX$ be the induced automorphism of $\cX$, so $S$ is again invariant for $f_\cX$, which we write again as $f$.  Let us write the various intersection products with $S$ as: $\mu_1=S\cdot E_1$, $\mu_0=S\cdot\Sigma_0$, $\mu_2=S\cdot E_2$, $\mu_P=S\cdot P$, $\mu_\gamma=S\cdot\Sigma_\gamma$, $\mu_Q=S\cdot Q$.  Since $e_1,e_2\in\Sigma_0$, we have
$$\mu_1+\mu_0+\mu_2=m.$$
Now we also have
$$f:\ \ \Sigma_\beta\to E_2\to \Sigma_0\to E_1$$
so $\mu_\beta=\mu_2=\mu_0=\mu_1=\mu$ for some positive integer $\mu$, and $m=3\mu$.  Similarly, $p,e_2\in\Sigma_\beta$, so we conclude that $\mu_P+\mu_\beta+\mu_2=m$, and thus $\mu_P=\mu$.  Following the backward orbit of $P$, we deduce that $S\cdot f^jQ=\mu$ for all $0\le j\le n$.

Now recall that if $L\in H^{1,1}(\bP^2,\bZ)$ is the class of a line, then the canonical class of $\bP^2$ is $-3L$.  Thus the canonical class $K_\cX$ of $\cX$ is $-3L+\sum E$, where sum is taken over all blowup fibers $E$.  In particular, the class of $S$ in $H^{1,1}(\cX)$ is $-\mu K_\cX$.  Since we obtained $\cX$ by performing $n+3$ blowups on $\bP^2$, the (arithmetic) genus formula, applied to the strict transform of $S$ inside $\cX$, gives:
$$g(S) = {S\cdot(S+K_\cX)\over 2} + 1 ={\mu(\mu-1)\over 2}K^2_\cX+1={\mu(\mu-1)\over 2}(9-(3+n))+1.$$
Now let $\nu$ denote the number of connected components of $S$ (the strict transform inside $\cX$).  We must have $g(S)\ge 1-\nu$.  Further, the degree $3\mu$ of $S$ must be at least as large as $\nu$, which means that $\mu(\mu-1)(n-6)\le 2\nu\le 6\mu$ and therefore $\mu\le 6/(n-6)+1$.  We have two possibilities: (i) If $n\ge 13$, then $\mu=1$, and $S$ must have degree 3; (ii) if $n=11$ or 12, either $\mu=1$, (i.e., the degree of $S$ is 3), or $\mu=2$.  Let us suppose $n=11$ or 12 and $\mu=2$.  From the genus formula we find that $5\le \nu\le 6$.  We treat these two cases separately.

{Case 1.}  $S$ cannot have 6 connected components.  Suppose, to the contrary, $\nu=6$.  First we claim that $S$ must be minimal, that is, we cannot have a nontrivial decomposition$S=S_1\cup S_2$, where $S_1$ and $S_2$ are invariant.   By the argument above,  $S_1$ and $S_2$ must be cubics, and thus they must both  contain all $n+3\ge 14$ points of blowup.   But then they must have a common component, so $S$ must be minimal.  

Since the degree of $S$ is 6, it follows that $S$ is the union of 6 lines which map $L_1\to L_2\to\cdots \to L_6\to L_1$.  Further, each $L_i$ must contain exactly one point of indeterminacy, since it maps forward to a line and not a quadric.  Since the class of $S$ in $H^{1,1}(\cX)$ is $-2K_\cX$, we see that $e_1,e_2,p,q\in S$ with multiplicity 2.  Without loss of generality, we may assume that $e_1\in L_1$, which means $\Sigma_\beta\cap L_1\ne\emptyset$, and therefore $e_2\in L_2$.  Similarly $q\in L_3$.  Since the backward image of $L_3$ is a line,  $e_1,e_2\notin L_3$, and thus $p\in L_3$, which gives $\Sigma\cap L_3\ne\emptyset$.  Continuing this procedure, we end up with $L_6\ni p,q$.  It follows that $L_1=L_4$, $L_2=L_5$, and $L_3=L_6$, so $S$ has only 3 components.

{Case 2.}  $S$ cannot have 5 connected components.  Suppose, to the contrary, that $\nu=5$.  It follows that $S$ is a union of 4 lines and one quadric.  Without loss of generality we may assume that $L_1\to Q\to L_2\to L_3$.  Since $L_1$ maps to a quadric, it cannot contain a point of indeterminacy, which means that $L_1\cap\Sigma_0\ne\emptyset$, $L_1\cap\Sigma_\beta\ne\emptyset$, and $L_1\cap\Sigma_\gamma\ne\emptyset$.  It follows that $e_1,e_2,q\in Q$.  On the other hand, since $L_2$ maps to a quadric by $f^{-1}$, we have $e_1,e_2,p\in Q$, and $e_1,e_2,q\notin L_2$.  Thus we have that $Q\cap \Sigma_0=\{e_1,e_2\}$, $Q\cap\Sigma_\beta=\{e_2,p\}$, and $Q\cap\Sigma_\gamma=\{e_1,p\}$.  It follows that $q\notin L_2$, which means that $L_2$ does not contain any point of indeterminacy, and therefore $L_2$ maps to a quadric.

Thus we conclude that $S$ has degree 3, so we may write $S=\{h=0\}$ for some cubic $h$.  Since the class of $S$ in $H^{1,1}(\cX)$ is $-3K_\cX$, we see that $e_1,e_2,q\in S$.  Since these are the images of the exceptional lines, the polynomial $h\circ f$ must vanish on $\Sigma_0\cup\Sigma_\gamma\cup\Sigma_\beta$.  Thus $j_f$ divides $h\circ f$, and since $h\circ f$ has degree 6, we must have (1.1).  \qed

\noindent{\bf Remarks. }  (a)  From the proof of Theorem 1.2, we see that  if $S$ is an invariant curve, $n\ge11$, then $S$ contains $e_1$, $e_2$, and $f^jq$, $0\le j\le n$.
(b) The only positive entropy parameters which are not covered in Theorem 1.2 are the cases $n=7,8,9,10$.  By Proposition B.1, we have $\cV_n\subset\Gamma$ for $7\le n\le 10$. 

\proclaim Corollary 1.3.  If $S$ is $f$-periodic with period $k$, and if $n\ge11$, then  $S\cup\cdots\cup f^{k-1}S$ is  invariant and thus a cubic.

\noindent{\bf \S2.  Invariant Cubics.  }   In this section, we identify the parameters $(a,b)\in\cV$ for which the birational map $f_{a,b}$ has an invariant curve, and we look at the behavior of $f_{a,b}$ on this curve.
We define the functions:
$$\eqalign{&\varphi_1(t)=\left({t-t^3-t^4\over 1+2t+t^2}, {1-t^5\over t^2+t^3}\right),\cr
\varphi_2(t)=&\left({t+t^2+t^3\over 1+2t+t^2}, {-1+t^3\over t+t^2}\right),\ \ \ \varphi_3(t)=\left(1+t,t- t^{-1}\right),\cr}\eqno(2.1)$$
The proofs of the results in this section involve some calculations that are possible but tedious to do by hand; but they are not hard with the help of Mathematica or Maple. 
\proclaim Theorem 2.1.   Let $t\ne0,\pm1$ with $t^3\ne1$ be given.  Then there is a homogeneous cubic polynomial $P$ satisfying (1.1) if and only if  $(a,b)=\varphi_j(t)$ for some $1\le j\le 3$.   If this occurs, then (up to a constant multiple) $P$ is given by (2.2) below.

\noindent{\it Proof. }  From the proof of Theorem 1.2, we know that $P$ must vanish at $e_1,e_2,q,p$.  Using the conditions $P(e_1)=P(e_2)=P(q)=0$, we may set
$$\eqalign{P[x_0:x_1:x_2] =& (-a^2C_1+aC_2)x_0^3 + C_2 x_1x_0^2+C_3x_2x_0^2+C_1x_0x_1^2+\cr
& +C_4x_2x_1^2+C_5x_0x_2^2+C_6 x_1x_2^2+C_7x_0x_1x_2\cr}$$
for some $C_1,\dots,C_7\in\bC$.  Since $e_1,e_2,q\in\{P=0\}$, we have $P\circ f=j_f\cdot \tilde P$ for some cubic $\tilde P$.  A computation shows that
$$\eqalign{ \tilde P= & (-ab^2C_1+b^2C_2+bC_3+aC_5)x_0^3 + (-2abC_1+2bC_2+C_3)x_1x_0^2+\cr
& (bC_1+C_5+aC_6+bC_7)x_2x_0^2 + (-aC_1+C_2)x_0x_1^2+C_1x_2x_1^2+(bC_4+C_6)x_0x_2^2+\cr
&+C_4x_1x_2^2 +(2bC_1+C_7)x_0x_1x_2.   }$$
Now setting $\tilde P=t P$ and comparing coefficients, we get a system of 8 linear equations in $C_1,\dots,C_7$ of the form
$$M\cdot[x_0^3,x_1x_0^2,x_2x_0^2,x_0x_1^2,x_0x_2^2,x_0x_1x_2,x_1x_2^2,x_1^2x_2]^t=0.$$
We check that there exist cubic polynomials satisfying (1.1) if and only if the two principal minors of $M$ vanish simultaneously, which means that
$$b(a+abt+abt^4-b^2t^4 -at^5+bt^5)=0$$
$$ -1+(1-a-b)t+(a+b)t^2+b^2t^3+b^2t^4+(a-2b)t^5+(1-a+2b)t^6-t^7=0 $$
Solving these two equations for $a$ and $b$, we obtain $\varphi_j$, $j=1,2,3$ as the only solutions, and then solving $M=0$ we find that $P$ must have the form:
$$\eqalign{P_{t,a,b}(x)=&a x_0^3 (-1 + t) t^4 +   x_1 x_2 (-1 + t) t (x_2 +  x_1 t) \cr
&+ x_0 [2 b  x_1 x_2 t^3 + 
              x_1^2 (-1 + t) t^3 +  x_2^2 (-1 + t) (1 + b t)]\cr
&+ x_0^2(-1 + t) t^3 [ a (x_1 + x_2 t) + t(x_1 +(-2b + t)x_2) ]   }\eqno(2.2)$$
which completes the proof. \qed

\noindent{\bf Remark.}  Let us discuss the values of $t$  omitted in Theorem 2.1.  There is no nonzero solution to (1.1) if $t=0$.  If $t=1$ or $t=-1$, then (1.1) is solvable iff $(a,b)\in\{b=0\}\subset\cV_0\cup\cV_1\cup\cV_6$.  If $t=\omega$ is a primitive cube root of unity, then $(a,b)=(0,0)=\varphi_2(\omega)$ or $(a,b)=(1+\omega,\omega-\bar\omega)=\varphi_3(\omega)$.  By [BK], $\varphi_3(\omega)\notin\cV_0\cup\dots\cup\cV_6$, and by Theorem 3.1, $\varphi_3(\omega)\notin\bigcup_{n\ge7}\cV_n$.

\medskip

If $h$ satisfies (1.1), then we may define a meromorphic 2-form $\eta_P$ on ${\bf P}^2$ by setting  $\eta_h:= {dx\wedge dy\over h(1,x,y)}$ on the affine coordinate chart $[1:x:y]$.  Then $\eta_h$ satisfies $t\, f^*\eta_h=\eta_h$.   It follows that if the points $\{p_1,\dots,p_k\}$ form a $k$-cycle which is disjoint from $\{h=0\}$, then the Jacobian determinant of $f$ around this cycle will be $t^{-k}$. 

Let  $\Gamma_j=\{(a,b)=\varphi_j(t): t\in\bC\}\subset\cV$ denote the curve corresponding to $\varphi_j$, and set $\Gamma:=\Gamma_1\cup\Gamma_2\cup\Gamma_3$.   Consistent with  [DJS], we find that the cases $\Gamma_j$ yield cubics with cusps, lines tangent to quadrics, and three lines passing through a point.

\epsfysize=2.0in
\centerline{ \epsfbox{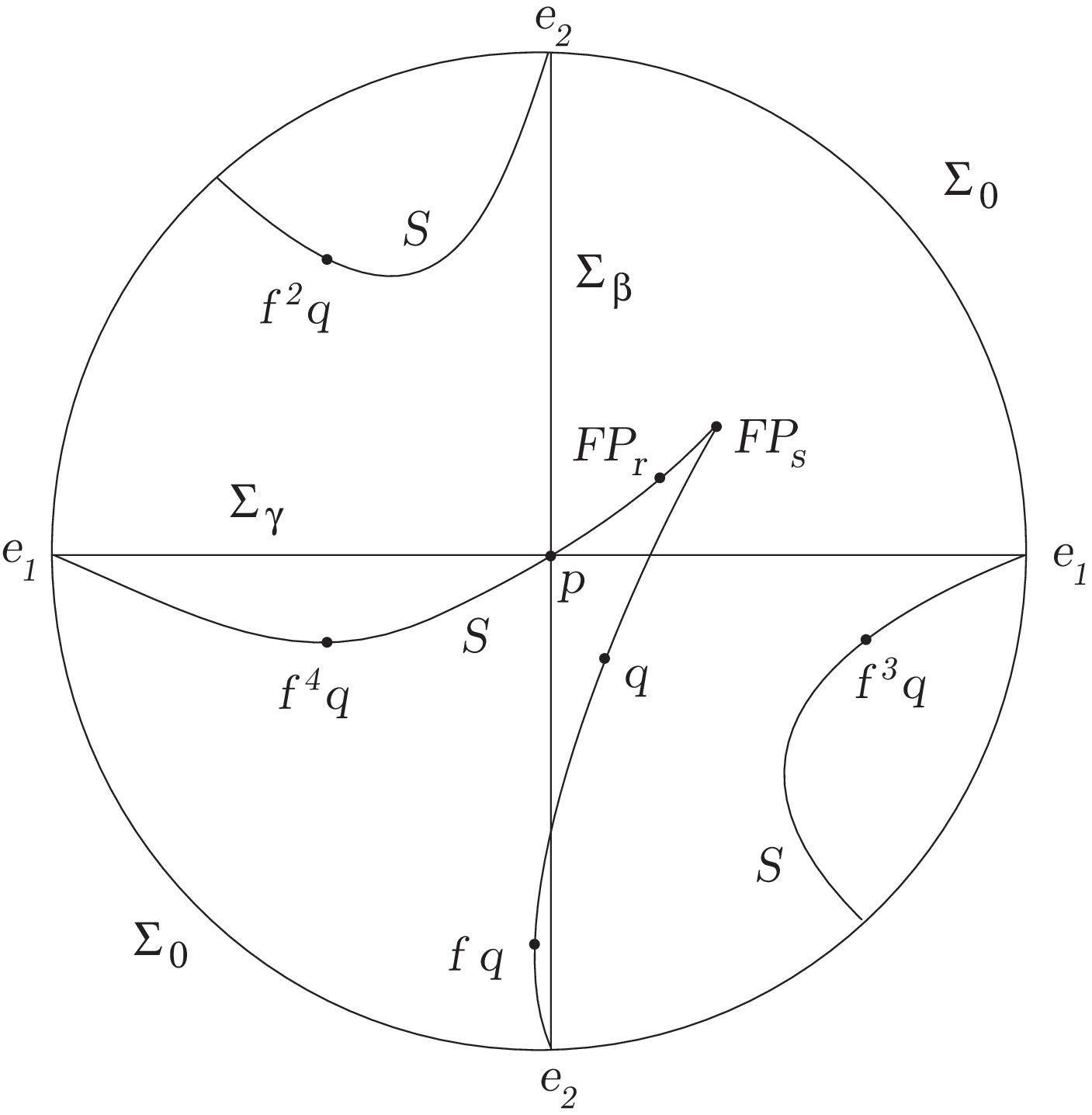}  }

\centerline{Figure 2.1.  Orbit of $q$ for family $\Gamma_1$; $1<t<\delta_\star$.}

\medskip\noindent{\it $\Gamma_1$: Irreducible cubic with a cusp. }  To discuss the family $\Gamma_1$, let $(a,b)=\varphi_1(t)$ for some $t\in\bC$.  Then the fixed points of $f_{a,b}$ are $FP_s=(x_s,y_s)$, $x_s=y_s=t^3/(1+t)$ and $FP_r=(x_r,y_r)$, $x_r=y_r=(-1+t^2+t^3)/(t^2+t^3)$.  The eigenvalues of $Df_{a,b}(FP_s)$ are $\{t^2,t^3\}$.    The invariant curve is $S=\{P_{t,a,b}=0\}$, with $P$ as in (2.2).   This curve $S$ contains $FP_s$ and $FP_r$, and has a cusp at $FP_s$.   The point $q$ belongs to $S$, and thus the orbit $f^jq$ for all $j$ until possibly we have $f^j q\in\cI$. The 2-cycle and 3-cycle are disjoint from $S$, so the multipliers in (B.1) must satisfy $\mu_2^3=\mu_3^2$, from which we determine that $\Gamma_1\subset\cV$ is a curve of degree 6.

We use the notation $\delta_\star$ for the real root of $t^3-t-1$.  Thus $1\le\lambda_n<\delta_\star$, and the $\lambda_n$ increase to $\delta_\star$ as $n\to\infty$.  The intersection of the cubic curve with $\bR\bP^2$ is shown in Figure 2.1.  The exceptional curves $\Sigma_\beta$ and $\Sigma_\gamma$ are used as axes, and we have chosen a modification of polar coordinates so that $\Sigma_0$, the line at infinity, appears as the bounding circle of $\bR\bP^2$.  The points $FP_{s/r}, e_1,e_2,p,q,fq,f^2q,f^3q$ all belong to $S$, and Figure 2.1 gives their relative positions with respect to the triangle $\Sigma_\beta,\Sigma_\gamma,\Sigma_0$  for all $1<t<\delta_\star$.  Since $t>1$, the points $f^jq$ for $j\ge4$ lie on the arc connecting $f^4q$ and $FP_r$, and $f^jq$ approaches $FP_r$ monotonically along this arc as $j\to\infty$.  In case  $(a,b)$ belongs to $\cV_n$, then $f^nq$ lands on $p$.  The relative position of $S_t$ with respect to the axes is stable for $t$ in a large neighborhood of $[1,\delta_\star]$.   However, as $t$ increases to $\delta_\star$, the fixed point $FP_r$ moves down to $p$; and for $t>\delta_\star$, $FP_r$ is in the third quadrant.   And as $t$ decreases to $1$,  $f$ approaches the (integrable) map  $(a,b)=(-1/4,0)\in\cV_6$.  The family $\cV_6$ will be discussed in Appendix A.  When $0<t<1$, the point $FP_s$ becomes attracting, and the relative positions of $q$ and  $fq$, etc., are reversed.  Figure 2.1 will be useful in explaining the graph shown in Figure 6.1.

\medskip\noindent{\it $\Gamma_2$: Line tangent to a quadric. }  Next we suppose that $(a,b)=\varphi_2(t)$.  We let $S=\{P_{t,a,b}=0\}$ be the curve in (2.2).  In this case, the curve is the union of a line $L=\{t^2x_0+tx_1+x_2=0\}$ and a quadric $Q$.  The fixed points are $FP_s=(x_s,y_s)$, $x_s=y_s=-t^2/(1+t)$ and $FP_r=(x_r,y_r)$, $x_r,y_r=(1+t+t^2)/(t+t^2)$.  The eigenvalues of $Df_{a,b}$ at $FP_s$ are $\{-t,-t^2\}$.  The 3-cycle and $FP_r$ are disjoint from $S$, so we have $\det(Df_{a,b}FP_r)^3-\mu_3=0$ on $\Gamma_2$, with $\mu_3$ as in (B.1).  Extracting an irreducible factor, we find that $\Gamma_2\subset\cV$ is a quartic.  

\epsfysize=2.0in

\centerline{ \epsfbox{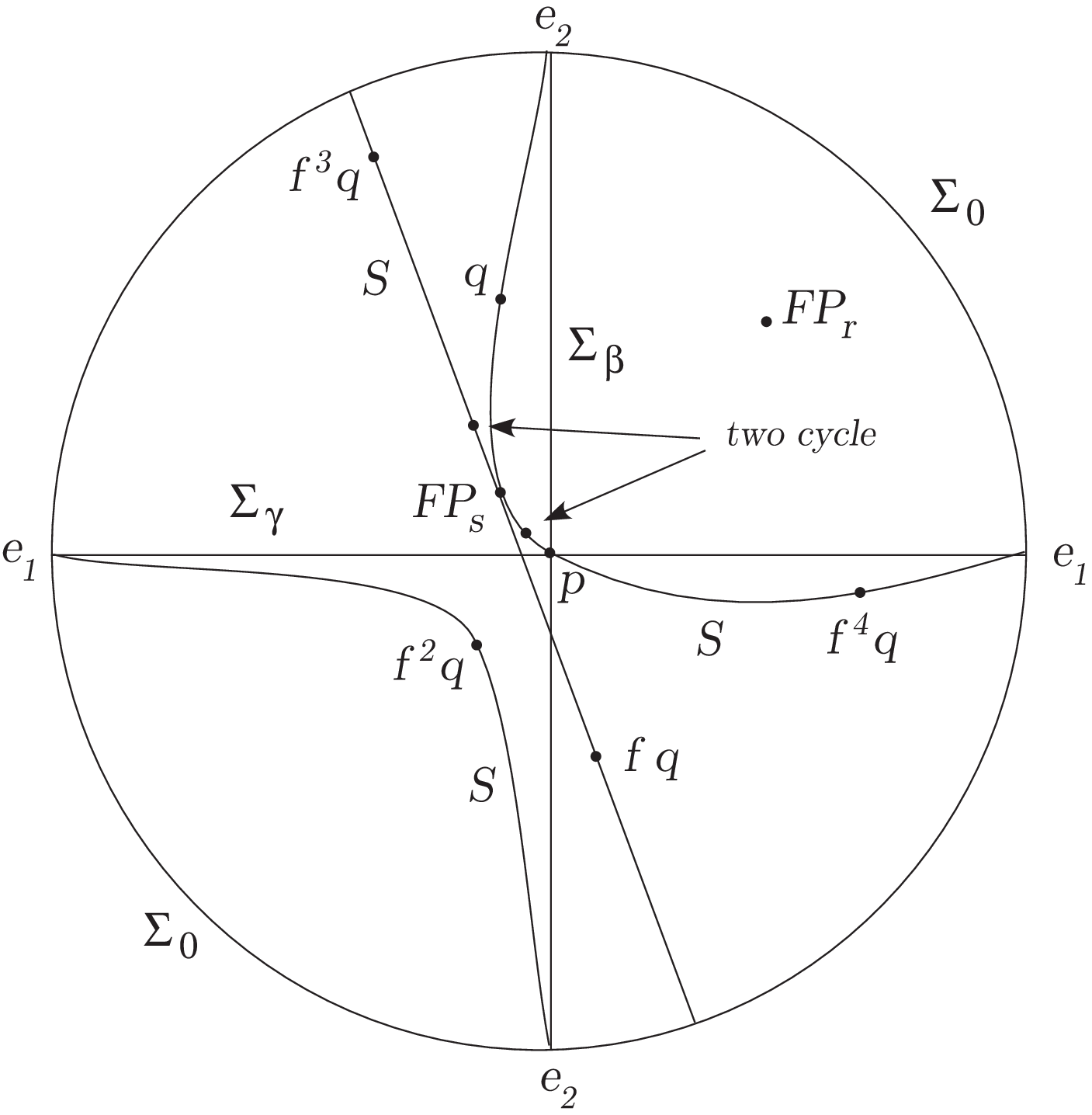}  }

\centerline{Figure 2.2.  Orbit of $q$ for family $\Gamma_2$; $1<t<\delta_\star$.}
\smallskip

Figure 2.2 gives for $\Gamma_2$ the information analogous to Figure 2.1.  The principal difference with Figure 2.1 is that $S$ contains an attracting 2-cycle; there is a segment $\sigma$  inside the line connecting $f^3q$ to one of the period-2 points, and there is  an arc $\gamma\ni p$ inside the quadric connecting $f^4q$ to the other period-2 point.  Thus the points $f^{2j+1}q$ will approach the two-cycle monotonically inside $\sigma$ as $j\to\infty$, and the points $f^{2j}q$ will approach the two-cycle monotonically inside $\gamma$.  The picture of $S$ with respect to the triangle $\Sigma_\beta,\Sigma_\gamma,\Sigma_0$ is stable for $t$ in a large neighborhood of $[1,\delta_\star]$.  As $t$ increases to $\delta_\star$, one of the points of the 2-cycle moves down to $p$.  As $t$ decreases to 1, $q$ moves up (and $f^2q$ moves down) to $e_2\in\cI$, and $fq$ moves down to $\Sigma_0$.  The case $t=1$ is discussed in Appendix A.

\epsfysize=2.0in
\centerline{ \epsfbox{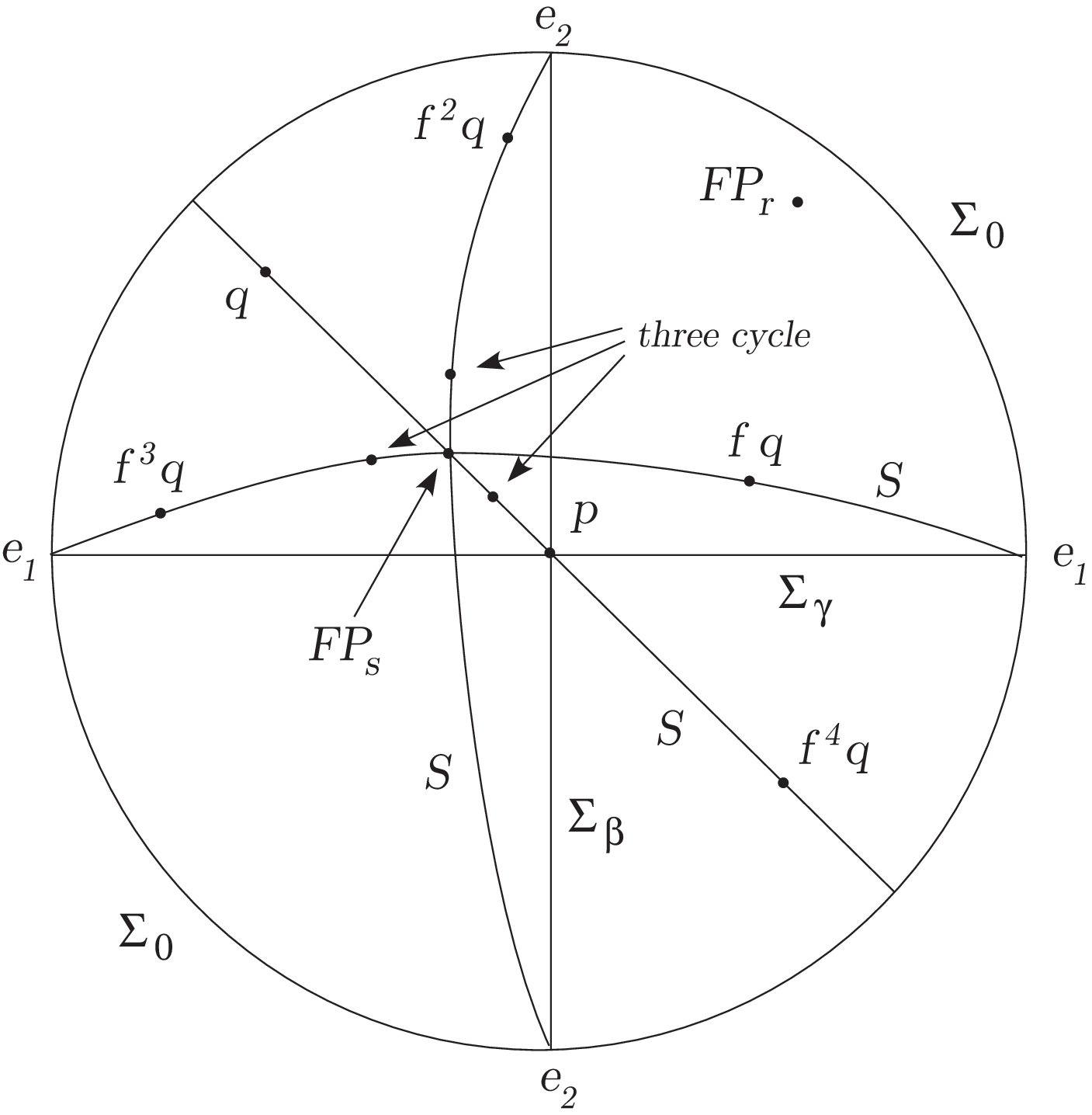}  }

\centerline{Figure 2.3.  Orbit of $q$ for family $\Gamma_3$; $1<t<\delta_\star$.}

\medskip\noindent {\it $\Gamma_3$: Three lines passing through a point. }  Finally, set $(a,b)=\varphi_3(t)$, and let 
$S=\{P_{t,a,b}=0\}$
be given as in (2.2).  The fixed points are $FP_s=(x_s,y_s)$, $x_s=y_s=-t$ and $FP_r=(x_r,y_r)$, $x_r=y_r=1+t^{-1}$.  The invariant set $S$ is the union of three lines $L_1=\{tx_0+x_1=0\}$, $L_2=\{tx_0+x_2=0\}$, $L_3=\{(t+t^2)x_0+tx_1+x_2=0\}$,  all of which pass through $FP_s$.   Further $p,q\in L_3\to L_2\to L_1$.  The eigenvalues of $Df_{a,b}$ at $FP_s$ are $\{\omega t,\omega^2 t\}$, where $\omega$ is a primitive cube root of unity.  The 2-cycle and $FP_r$ are disjoint from $S$, so we have ${\rm det}(Df_{a,b}FP_r)^2-\mu_2=0$ on $\Gamma_2$.  Extracting an irreducible factor from this equation we see that $\Gamma_3\subset\cV$ is a quadric.  Figure 2.3 is analogous to Figures 2.1 and 2.2; the lines $L_1$ and $L_2$ appear curved because of the choice of coordinate system.

\bigskip\noindent{\bf \S3.  Roots and Factorization of $\chi_n$. }  There is a close connection (Theorem 3.1) between the parameters $(a,b)=\psi_j(t)$ and the roots $t$ of $\chi_n$.   In fact we may write $\chi_n=C_n\psi_n$, where $C_n$ is a product of cyclotomic factors, and $\psi_n$ is the minimal polynomial of $\lambda_n$.  In Theorems 3.3 and 3.5 we give the precise form of the cyclotomic part $C_n$ of this factorization.   One consequence of the factorization of $\chi_n$ is to obtain a precise count of the number of automorphisms in ${\cal V}_n$ with invariant curves.

\proclaim Theorem 3.1.  Suppose that $n$, $1\le j\le 3$, and $t$ are given, and suppose that $(a,b):=\varphi_j(t)\notin\cV_k$ for any $k<n$.  Then the point $(a,b)$ belongs to $\cV_n$ if and only if: $j$ divides $n$ and $t$ is a root of $\chi_n$.

\noindent{\it Proof. }  Let us start with the case $j=3$ and set $(a,b)=\varphi_3(t)$.  By the calculations in \S2, we know that $S=L_1\cup L_2\cup L_3$ factors into the product of lines, each of which is invariant under $f^3$.  $L_3$ contains $FP_s$ and $R=[t^2:-1:t-t^3-t^4]$, which is periodic of period 3.  We define $\psi(\zeta)=FP_s+\zeta R$, which gives a parametrization of $L_3$; and the points $\psi(0)$ and $\psi(\infty)$ are fixed under $f^3$.  The differential of $f^3$ at $FP_s$  was seen to be $t^3$ times the identity, so we have $f^3(\psi(\zeta))=\psi(t^3\zeta)$.  Now set $\zeta_q:=t^2/(1-t^2-t^3)$ and $\zeta_p:= t/(t^3-t-1)$.  It follows that $\psi(\zeta_q)=q$ and $\psi(\zeta_p)=p$.   If $n=3k$, then $f^nq=f^{3k}q=p$ can hold if and only if $t^n\zeta_q=t^{3k}\zeta_q=\zeta_p$, or 
$t^{n+2}/(1-t^2-t^3)=t/(t^3-t-1)$,  which is equivalent to $\chi_{n}(t)=0$.

Next, suppose that $j=2$ and let $(a,b)=\varphi_2(t)$.  In this case the polynomial $P$ given in (3.1) factors into the product of a line $L$ and a quadric $Q$.  $L$ contains $FP_s$ and the point $R=[t+t^2:t^3+t^2-1:-t]$, which has period 2.  We parametrize $L$ by the map $\psi(\zeta)=FP_s+\zeta R$.  Now $f^2$ fixes $FP_s$ and $R$, and the differential of $f^2$ has an eigenvalue $t^2$ in the eigenvector $L$, so we have $f^2\psi(\zeta)=\psi(t^2\zeta)$.  Since $p,q\in Q$, we have $fq,f^{-1}p\in L$.  We see that $\zeta_q:=t^3/(1-t^2-t^3)$ and $\zeta_p:=(t^3-t-1)^{-1}$ satisfy $\psi(\zeta_q)=fq$ and $\psi(\zeta_p)=f^{-1}p$.  If $n=2k$, then the condition $f^nq=f^{2k}q=p$ is equivalent to the condition $t^{2n-2}\zeta_q=\zeta_p$, which is equivalent to $\chi_n(t)=0$.

Finally we consider the case $j=1$ and set $(a,b)=\varphi_1(t)$.   If we substitute these values of $(a,b)$ into the formula (3.1), we obtain a polynomial $P(x)$ which is cubic in $x$ and which has coefficients which are rational in $t$.  In order to parametrize $S$ by $\bC$, we set $\psi(\zeta)=FP_s+\zeta A+\zeta^2 B +\zeta^3 FP_r$.  We may solve for $A=A(t)$ and $B=B(t)$ such that $P(\psi(\zeta))=0$ for all $\zeta$.   Thus $f$ fixes $\psi(0)$ and $\psi(\infty)$, and $f(\psi(\zeta))=\psi(t\zeta)$.  We set $\zeta_q:=t^2/(1-t^2-t^3)$ and $\zeta_p:=t/(t^3-t-1)$.  The condition $f^nq=p$ is equivalent to $t^n\zeta_q=\zeta_p$, or $-t^{n+2}/(t^3+t^2-1)=t/(t^3-t-1)$, or $\chi_n(t)=0$.  \qed

For each $n$, we let $\psi_n(t)$ denote the minimal polynomial of $\lambda_n$.   
\proclaim Theorem 3.2.  Let $t\ne1$ be a root of $\chi_n$ for $n\ge7$.  Then either $t$ is a root of $\psi_n$, or $t$ is a root of $\chi_j$ for some $0\le j\le 5$.

\noindent{\it Proof. }  Let $t$ be a root of $\chi_n$ which is not a root of $\psi_n$.  This root has modulus 1, so by Lemma 1.6 of [W], it must be a root of unity.  We will show it is a root of $\chi_j$ for some $0\le j\le 5$.  First we note that $\chi_6(t) = (t-1)^3(t+1)(t^2+t+1)(t^4+t^3+t^2+t+1)$, and so every root of $\chi_6$ is a root of $\chi_j$ for some $0\le j\le 5$.  Similarly, we may check that the Theorem is true for $7\le n\le 13$.  By induction, it suffices to show that if $t$ is a root of unity, then it is a root of $\chi_j$ for some $0\le j\le n-1$.  By Theorem 3.1, we see that $\chi_n(t)=0$ if and only if $t^n\zeta_q(t)=\zeta_p(t)$, where $\zeta_q(t)=t^2/(1-t^2-t^3)$ and $\zeta_p(t)=t/(t^3-t-1)$. 

Claim 1: We may assume $t^n\ne\pm1$.  Otherwise, from $t^n\zeta_q(t)=\zeta_p(t)$ we have 
$$t^2(t^3-t-1)\pm t(t^3+t^2-1)=0.$$
 In  case we take ``+'', the roots are also roots of $\chi_0$, and in  case we take ``$-$'', there are no roots of unity.

Claim 2: If $t^k=1$ for some $0\le k\le n-1$, then $t$ is a root of $\chi_j$ for some $0\le k-1$.  As in the proof of Theorem 2.2, the orbit of $\zeta_q$ is $\{\zeta_q,t\zeta_q,\dots,t^{k-1}\zeta_q\}$.  Thus the condition that $t^k=1$ means that $f^kq=p$, so $\chi_k(t)=0$.

Claim 3:   If $t$ is a primitive $k$th root of unity, we cannot have $k>n$.  Otherwise, we have $t^m\zeta_p=\zeta_q$ for $m=k-n$, which means that $t$ is a root of $\xi_m:=t^{m-1}(t^3+t^2-1)+t^3-t-1$.  This means that $t$ is a root of both polynomials:
$$\eqalign{\chi_n/(t-1) = &1+t-(t^3+t^4+\cdots+t^n) + t^{n+2}+t^{n+3}\cr
\xi_m/(t-1) = &1 + 2t + 2t^2 + (t^3+\cdots + t^{n-2}) + 2t^{m-1} + 2t^m + t^{m+1}.}$$

We divide our argument into four sub-cases: $1\le m\le 4$, $5\le m\le n-2$, $n-1\le m\le n+6$, and $n+6<m$.  If $1\le m\le 4$, we may check directly that $\xi_m$ does not have  a $k$th root of unity with $k>n$.  If $5\le m\le n-2$, then ${k}={n+m} \ge {2m+2}$.  We may assume that $t=e^{ i\theta}$, with $\theta=2\pi/k$.  Thus $\Im m(t^j)>0$ for $j=1,\dots,m+1$.  It follows that $\Im m(\xi_m(t)/(t-1))>0$, contradicting the assumption that $t$ is a root.  

If $n-1\le m\le n+6$, we show that $\xi_m$ does not vanish at a $k$-th root of unity by inspecting the remainder obtained after performing three steps of the Euclidean Algorithm on the polynomials $\chi_n$ and $\xi_m$.

This leaves the case  $2n+6<k$, i.e.\ $m>n+6$.  Here we have
$$\eqalign{\xi_m/(t-1) + \chi_n/&(t-1) = 2 + 3t + 2t^2 + t^{n+1}+ 2t^{n+2} + 2t^{n+3} +\cr
&+ t^{n+4} + \cdots+t^{m-4} + t^{m-2} + 2t^{m-1} + 2t^m + t^{m+1}.}$$
Since $t^k=1$, and since $m=n-k$, we have $\Im m(t^{n+1}+t^{m-1}) = \Im m(t^{n+2} + t^{m-2})=\cdots=0$.  Applying this repeatedly, we find
$$\eqalign{\Im m(\xi_m/(t-1) + \chi_n/(t-1))= \Im m& (3t+2t^2+t^{n+2} + t^{n+3} + t^{m-1}+ 2t^m+t^{m+1})\cr
=\Im m& (3t + 2t^2 -t^{n-1}-2t^n-t^{n+1} + t^{n+2} + t^{n+3})=0.
}$$
Again taking $t=e^{ i\theta}$ with $\theta=2\pi /k$, the previous equation becomes
$$\eqalign{3\sin\theta +  2&\sin(2\theta) -\sin((n-1)\theta) -2\sin (n\theta)-\sin((n+1)\theta)+\cr
& +\sin((n+2)\theta) + \sin((n+3)\theta) \cr
 =   3 &\sin\theta +   2\sin(2\theta) + 2\left(\cos((n+1)\theta)  + \cos((n+2)\theta)\right)\sin\theta \cr
&  -\sin((n-1)\theta) -\sin (n\theta)  .}$$
By  $2\theta/\pi<\sin\theta=\sin(\pi-\theta)<\theta$ for $0<\theta<\pi/2$, this last expression is seen to be strictly negative since $\theta=2\pi/k$ with  $k>2n+6$ and $n\ge13$.  This completes the proof of Claim~3.
 \qed
\noindent{\bf Remark. }  Let us recall the construction of surface automorphisms given in [M2].  It starts with a root $t$ of $\psi_n$ and a singular cubic, which is one of the three cases considered above: cusp cubic, line and quadric, or three lines.  The automorphism is then uniquely determined by a certain marking of the cubic.  Now if we choose $(a,b)=\psi_*(t)$, then the centers of blowup on the invariant curve in $\cX_{a,b}$ correspond exactly to the marked points, and  $f_{a,b}$ is equivalent to the resulting automorphism.
\medskip
The following result gives the possibilities for the roots of $\chi_n(x)/\psi_n(x)\in\bZ[x]$.
\proclaim Theorem 3.3.  Let $t\ne1$ be a root of $\chi_n$ with $n\ge7$.  Then $t$ is either a root of $\psi_n$, or $t$ is a root of some $\chi_j$ for $0\le j\le 5$.  Specifically, if $t$ is not a root of $\psi_n$, then it is a $k$th root of unity corresponding to one of the following possibilities:
\item{(i)}  $k=2$, $t+1=0$, in which case 2 divides $n$;
\item{(ii)} $k=3$,  $t^2+t+1=0$, in which case 3 divides $n$;
\item{(iii)} $k=5$, $t^4+t^3+t^2+t+1=0$, in which case $n\equiv 1{\rm \ mod\ }5$;
\item{(iv)} $k=8$, $t^4+1=0$, in which case $n\equiv2{\rm \ mod \ }8$;
\item{(v)} $k=12$,  $t^4-t^2+1=0$, in which case $n\equiv 3{\rm \ mod\ }12$;
\item{(vi)} $k=18$, $t^6-t^3+1=0$,  in which case $n\equiv 4{\rm \ mod\ }18$;
\item{(vii)} $k=30$, $t^8+t^7-t^5-t^4-t^3+t+1=0$,  in which case $n\equiv 5{\rm \ mod\ }30$.
\vskip0pt \noindent Conversely,  for each $n\ge7$ and $k$ satisfying one of the conditions above, there is a corresponding root $t$ of $\chi_n$ which is a $k$th root of unity.

\noindent{\it Proof. }  Recall that $\chi_n(t)=0$ if and only if $t^n\zeta_q(t)=\zeta_p(t)$ by Theorem 2.2.  If $t$ is a $k$-th root of unity, then $k<n$ and $\chi_j(t)=0$ for some $0\le j\le 5$.  In case $j=0$, $\zeta_p(t)=\zeta_q(t)$, and $(t+1)(t^2+t+1)=0$.  Thus $t^n\zeta_q(t)=\zeta_p(t)$ if and only if $t+1=0$ and 2 divides $n$, or $t^2+t+1=0$ and 3 divides $n$.  Now let us write $k_j=5,8,12,18,30$ for $j=1,2,3,4,5$, respectively, in case we have $1\le j\le 5$, $t^j\zeta_q(t)=\zeta_p(t)$, and $n\equiv j{\rm\ mod\ }k_j$.  Thus  $t^n\zeta_q(t)=\zeta_p(t)$ if and only if $n\equiv j{\rm\ mod\ }k_j$.  That is, 
$t^n\zeta_q(t)=(t^{k_j})^nt^j\zeta_q(t)=t^j\zeta_q(t)=\zeta_p(t)$.  \qed

As a corollary, we see that the number of elements of $\Gamma_j\cap\cV_n$ is determined by the number of Galois conjugates of $\lambda_n$.
\proclaim Corollary 3.4.  If $n\ge7$, and if $1\le j\le 3$ divides $n$, then 
$$\Gamma_j\cap\cV_n=\{\varphi_j(t): t {\rm\ is\ a\ root\ of\ }\psi_n\}.$$
In particular, these sets are nonempty.

\proclaim Theorem 3.5.  If $n\ge7$, then every root of $\chi_n$ is simple.  Thus the possibilities enumerated in Theorem 3.3 give the irreducible factorization of $\chi_n$.  In particular, since each of the factors in (i--vii) below can occur only once, we have $n-26\le {\rm deg}(\psi_n)\le n+3$.

\noindent{\it Proof. }  If $t$ is a root of $\chi_n$, then either it is a root of $\psi_n$, which is irreducible, or it is one of the roots of unity listed in Theorem 3.3.  We have 
$$\chi'_n(t)=(n+4)t^{n+3}-(n+2)t^{n+1}-(n+1)t^n+3t^2+2t.$$
Since $\chi'_n(1)=6-n$, 1 is a simple root.  Now we check all the remaining cases:
\item{(i)} 2 divides $n$ :  $t+1=0\Rightarrow\chi'_n(t)=-2-n\ne0$.
\item{(ii)} 3 divides $n$ :  $t^2+t+1=0\Rightarrow\chi'_n(t)=3t^2-nt+3\ne0$.
 \item{(iii)}  $n\equiv1$ (mod 5) : $t^4+t^3+t^2+t+1=0\Rightarrow \chi'_n(t)=(n+4)t^4-(n-1)t^2-(n-1)t\ne0$.
 \item{(iv)} $n\equiv 2$ (mod 8) : $t^4+1=0\Rightarrow \chi'_n(t)=-(n+2)t^3-(n-2)t^2-(n+2)t\ne0$.
\item{(v)} $n\equiv 3$ (mod 12) : $t^4-t^2+1=0\Rightarrow \chi'_n(t)=-(n+1)t^3-(n-1)t^2+2t-2\ne0$.
\item{(vi)} $n\equiv 4$ (mod 18) : $t^6-t^3+1=0\Rightarrow \chi'_n(t)=-(n+2) t^5+3t^4+3t^2-(n+2)t\ne0$. 
\item{(vii)} $n\equiv 5$ (mod 30) : $t^8+t^7-t^5-t^4-t^3+t+1=0\Rightarrow $
\vskip0pt \ \ \ \ \ \ \ $\chi'_n(t)=(n+4)t^8-(n+2)t^6-(n+1)t^5+3t^2+2t\ne0$.  \qed

\noindent{\bf Example. } The number $n=26$ corresponds to cases (i), (iii), and (iv), so we see that 
$\chi_{26}=(t-1)(t+1)(t^4+1)(t^4+t^3+t^2+t+1)\psi_{26},$
 so $\psi_{26} $
 has degree 20.

\bigskip\noindent{\bf\S4.  Surfaces without Anti-PluriCanonical Section. }   A curve $S$ is said to be a pluri-anticanonical curve if it is the zero set of a section of $\Gamma(\cX,(-K_\cX)^{\otimes n})$ for some $n>0$.  We will say that $(\cX,f)$ is minimal if whenever $\pi:\cX\to\cX'$ is a birational morphism mapping $(\cX,f)$ to an automorphism $(\cX',f')$, then $\pi$ is an isomorphism.  Gizatullin conjectured that if $\cX$ is a rational surface which has an automorphism $f$ such that $f^*$ has infinite order on $Pic(\cX)$, then $\cX$ should have an anti-canonical curve.  Harbourne [H] gave a counterexample to this, but this counterexample is not minimal, and $f$ has zero entropy.   See [Zh] for positive results in the zero entropy case.
\proclaim Proposition 4.1.  Let $\cX$ be a rational surface with an automorphism $f$.  Suppose that $\cX$ admits a pluri-anticanonical section.  Then there is an $f$-invariant curve.

\noindent{\it Proof. }  Suppose there is a pluri-anticanonical section.  Then $\Gamma(\cX,(-K_\cX)^{\otimes n})$ is a nontrivial finite dimensional vector space for some $n>0$, and $f$ induces a linear action on this space.  Let $\eta$ denote an eigenvalue of this action.  Since $\cX$ is a rational surface,  $S=\{\eta=0\}$ is a nontrivial curve, which must be invariant under $f$.  \qed
The following answers a question raised in [M2, \S12].
\proclaim Theorem 4.2.  There is a rational surface $\cX$ and an automorphism $f$ of $\cX$ with positive entropy such that $(\cX,f)$ is minimal, but there is no $f$-invariant curve.  In particular, there is no pluri-canonical section.

\noindent{\it Proof. }   We start by finding a map without invariant curve.  We consider $(a,b)\in \cV_{11}$.  Suppose that $(\cX_{a,b},f_{a,b})$ has an invariant curve $S$.  Then by Theorem 1.2, $S$ must be a cubic.  By Theorem 3.1  we must have $(a,b)\in\Gamma_1$.  That is, $(a,b)=\varphi_1(t)$ for some $t$.  By Theorem 3.3, $t$ is a root of the minimal polynomial $\psi_{11}$.  By Theorem 3.5, $\psi_{11}(t)=\chi_{11}(t)/((t-1)(t^4+t^3+t^2+t+1))$ has degree 10, so $\cV_{11}\cap\Gamma_1$ contains 10 elements.   However,  there are 12 elements in $\cV_{11}-\Gamma_1$.  Each of these gives an automorphism $(\cX_{a,b},f_{a,b})$ with entropy $\log\lambda_{11}>0$ and with no invariant curve.  

Next we claim that $(\cX_{a,b},f_{a,b})$ is minimal.  Suppose that $\pi:\cX_{a,b}\to\cY$ is a holomorphic map which is birational, and suppose that the induced map $f_\cY$ is an automorphism.  Then there are finitely many points $P=\{p_1,\dots,p_N\}$ such that $\pi^{-1}p_j$ has positive dimension.  It follows that $P$ is invariant under $f_\cY$ and thus $\pi^{-1}P$ is an invariant curve for $f_{a,b}$.  The nonexistence of invariant curves then shows that $(\cX_{a,b},f_{a,b})$ is minimal.  

The nonexistence of an pluri-anticanonical section follows from Proposition 4.1.  \qed

\noindent {\bf Remark.}  By Proposition B.1 we cannot take $n \le 10$ the proof of Theorem 4.2.

\bigskip\noindent{\bf \S5.  Rotation  Domains.  }  Given an automorphism $f$ of a compact surface $\cX$, we define the Fatou set $\cF$ to be the set of normality of the iterates $\{f^n:n\ge0\}$.  Let $\cD$ be an invariant  component of $\cF$.  We say that $\cD$ is a rotation domain if $f_{\cD}$ is conjugate to a linear rotation (cf.\ [BS1] and [FS], as well as Ueda [U], where this terminology is more completely justified).  In this case, the normal limits of $f^n|_{\cD}$ generate a compact abelian group.  In our case, the map $f$ does not have finite order, so the iterates generate a torus $\bT^d$, with $d=1$ or $d=2$.  We say that $d$ is the rank of $\cD$.  The rank is equal to the dimension of the closure of a generic orbit of a point of $\cD$.   Let us start by proving the existence of rank 1 rotation domains.  

\proclaim Theorem 5.1.  Suppose that $n\ge7$, $j$ divides $n$, and $(a,b)\in\Gamma_j\cap \cV_n$.  That is, $(a,b)=\varphi_j(t)$ for some $t\in \bC$.  If $t\ne\lambda_n,\lambda^{-1}_n$ is a Galois conjugate of $\lambda_n$, then $f_{a,b}$ has a rotation domain of rank 1 centered at $FP_s$.

\noindent{\it Proof.}  There are three cases.   We saw in \S1 that the eigenvalues of $Df_{a,b}$ at $FP_s$ are $\{t^2, t^3\}$ if $j=1$; they are $\{-t,-t^2\}$ if $j=2$ and $\{\omega t,\omega^2t\}$ if $j=3$.  Since $\lambda_n$ is a Salem number, the Galois conjugate $t$ has modulus 1.  Since $t$ is not a root of unity, it satisfies the Diophantine condition 
$$|1-t^k|\ge C_0k^{-\nu}\eqno(5.1)$$ 
for some $C_0,\nu>0$ and all $k\ge2$.  This is a classical result in number theory.  A more recent proof (of a more general result) is given in  Theorem 1 of  [B].  We claim now that if $\eta_1$ and $\eta_2$ are the eigenvalues of $Df_{a,b}$ at $FP_s$, then for each $m=1,2$, we have 
$$|\eta_m-\eta_1^{j_1}\eta_2^{j_2}|\ge C_0(j_1+j_2)^{-\nu}\eqno(5.2)$$
 for  some $C_0,\eta>0$ and all $j_1,j_2\in\bN$ with $j_1+j_2\ge2$.  There are three cases to check: $\Gamma_j$, $j=1,2,3$.  In case $j=1$, we have that $\eta_m-\eta_1^{j_1}\eta_2^{j_2}$ is equal to either $t^2-t^{2j_1+3j_2}=t^2(1-t^{2(j_1-1)+3j_2})$ or $t^3-t^{2j_1+3j_2}=t^3(1-t^{2j_1+3(j_2-1)})$.  Since $j_1+j_2>1$, we see that (5.2) is a consequence of (5.1).  In the case $j=2$, $Df^2$ has eigenvalues $\{t^2,t^4\}$, and in the case $j=3$, $Df^3$ has eigenvalues $\{t^3,t^3\}$.  In both of these cases we repeat the argument of the case $j=1$.  It then follows from Zehnder [Z2] that $f_{a,b}$ is holomorphically conjugate to the linear map $L={\rm diag}(\eta_1,\eta_2)$ in a neighborhood of $FP_s$.  \qed

\epsfysize=2.0in
\centerline{ \epsfbox{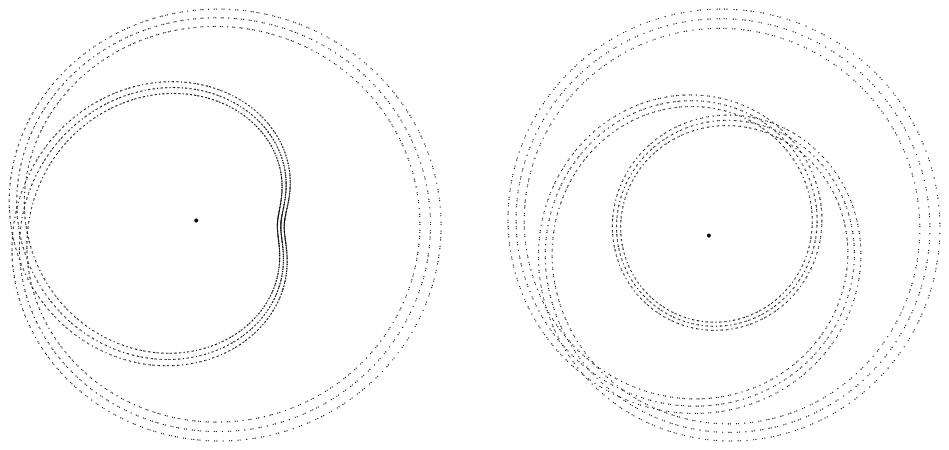}  }

\centerline{Figure 5.1.  Orbits of three points in the rank 1 rotation domain}

\centerline{ containing $FP_s$;  $(a,b)\in\cV_7\cap\Gamma_1$.  Two projections.}
\bigskip  For each $n\ge7$ and each divisor $1\le j\le 3$ of $n$, the only values of $(a,b)\in \cV_n\cap\Gamma_j$ to which Theorem 5.1 does not apply are the two values $\varphi_j(\lambda_n)$ and $\varphi_j(\lambda_n^{-1})$.  For all the other maps in $\cV_n\cap\Gamma$, the Siegel domain $\cD\ni FP_s$ is a component of both Fatou sets $\cF(f)$ and $\cF(f^{-1})$.  For instance, if $j=1$, then $f$ is conjugate on $\cD$ to the linear map $(z,w)\mapsto(t^2z,t^3w)$.  Thus, in the linearizing coordinate, the orbit of a point of $\cD$ will be dense in the curve $\{|z|=1,w^2=cz^3\}$, for some $r$ and $c$.  In particular, the closure of the orbit bounds an invariant (singular) complex disk.  Three such orbits are shown in Figure~5.1.  

   Now let us discuss the other fixed point.  Suppose that $(a,b)\in\Gamma_1\cap\cV_n$ and $\{\eta_1,\eta_2\}$ are the multipliers at $FP_r$.  As was noted in the proof above, $t$ satisfies (5.1), and so by Corollary B.5, both $\eta_1$ and $\eta_2$ satisfy (5.1).  On the other hand, the resonance given by Theorem B.3 means that they do not satisfy (5.2), and thus we cannot conclude directly that $f$ can be linearized in a neighborhood of $FP_r$.  However, by P\"oschel [P], there are holomorphic Siegel disks (of complex dimension one)  $s_j:\{|\zeta|<r\}\to \cX_{a,b}$, $j=1,2$, with the property that  $s_j'(0)$ is the $\eta_j$ eigenvector, and $f(s_j(\zeta))= s_j(\eta_j\zeta)$.  We note that one of these Siegel disks will lie in the invariant cubic itself and will thus intersect the rotation domain about $FP_s$.   By Theorem B.4 there are similar resonances between the multipliers for the 2- and 3-cycles, and thus similar Siegel disks, in the cases $\Gamma_2\cap\cV_n$ and $\Gamma_3\cap\cV_n$ respectively.

If  $j=2,3$ and $(a,b)\in\Gamma_j\cap\cV_n$, then at $FP_r$ we have
$$Df_{a,b}=\pmatrix{0&1\cr -{1\over t} & {1+t\over 1+t+t^2}\cr}{\rm \ for\ }j=2,{\rm\ and\ }=\pmatrix{0& 1\cr -{1\over t}  & {1\over t+1}\cr}\ {\rm \ for\ }j=3.\eqno(5.3)$$
We let $\eta_1,\eta_2$ denote the eigenvalues of $Df_{a,b}$ at $FP_r$.  For both $j=2$ and $j=3$ we have $\eta_1\eta_2={1/ t}$.
\proclaim Lemma 5.2.  Let $|t|=1$, let $(a,b)=\varphi_j(t)$, and let $\eta_1,\eta_2$ be the eigenvalues of $Df_{a,b}$ at $FP_r$.
\item{1.}  If $j=2$, then:\ \  $|\eta_1|=|\eta_2|=1\Leftrightarrow Re(t)\ge-7/8$.
\item{2.}  If $j=3$, then:\ \  $|\eta_1|=|\eta_2|=1\Leftrightarrow Re(t)\notin({-23-\sqrt{17}\over 32},{-23+\sqrt{17}\over 32})$.

\noindent{\it Proof. }  Set $t=e^{i\theta}$.  Since $\eta_1\eta_2=1/t$, we may set $\eta_1=\delta e^{i\omega}$, $\eta_2={1\over \delta}e^{-i(\theta+\omega)}$, $\delta>0$.  Consider first the case $j=3$. With $\eta_1,\eta_2$ written as above, it follows from (5.3) that they are the eigenvalues if and only if $\eta_1+\eta_2=1/(t+1)$.  Setting $\delta=1$, then, we have
$$\eqalign{\eta_1+\eta_2&={1\over t+1}\Leftrightarrow -1+2\cos \omega +2\cos(\theta+\omega)=0\cr
&\Leftrightarrow \cos\omega={1\over 4}\pm{\sqrt{ -(\cos^2\theta-1)(8\cos\theta+7)}\over 4(1+\cos\theta)}. \cr}$$
Thus there are solutions for $\omega$ if and only if the right hand side of the last equation is between $-1$ and $1$.  We may check that this happens if and only if $\cos\theta\ge-7/8$, which gives condition (2).

The case $j=2$ is similar.  Again we set $\delta=1$, and this time we find $\omega$ so that we have
$$\eqalign{&\eta_1+\eta_2={1+t\over 1+t+t^2} \Leftrightarrow -\cos{\theta\over2}+(1+2\cos\theta)\cos({\theta\over 2} +\omega)=0\cr
&\Leftrightarrow \cos\omega={\tau^2 \pm \sqrt{(1-\tau^2)(4\tau^2 + \tau-1)(4\tau^2 -\tau-1)}  \over 4\tau^2-1},\cr} $$
where we set $\tau=\cos{\theta\over 2}$ in the last equation.
Solving for the right hand side of this last equation to be between $-1$ and $1$, we find condition (1).  \qed

\proclaim Theorem 5.3.   Suppose that $n\ge8$, $j=2$ or 3, and $j$ divides $n$.  Suppose that $(a,b)=\varphi_j(t)$, and $|t|=1$ is a root of $\psi_n$ which satisfies the condition of Lemma 5.2.  Then $f_{a,b}$ has a rotation domain of rank 2 centered at $FP_r$.

\noindent{\it Proof. }  By Lemma 5.2, the multipliers $\eta_1,\eta_2$ of $Df_{a,b}$ at $FP_r$ both have modulus 1.  We will show that $\eta_1^{m_1}\eta_2^{m_2}\ne1$ for all $m_1,m_2\in\bZ$, not both zero.   Thus $\eta_1$ and $\eta_2$ are multiplicatively independent, and the rank of a rotation domain with such multipliers must be 2.   In particular, we then have the non-resonance condition, which asserts that  for $m=1,2$, we have  $\eta_m-\eta_1^{j_1}\eta_2^{j_2}\ne0$ for all $j_1,j_2\in\bN$ with $j_1+j_2\ge2$.  Since $\eta_1$ and $\eta_2$ are algebraic, it follows from Theorem 1 of [B] that (5.2) holds.  Thus by [Z2], $f_{a,b}$ is linearizable at $FP_r$.

To show multiplicative independence, let us start with the case $j=3$.   Suppose that $\eta_1^{m_1}\eta_2^{m_2}=1$ with $m_1\ge m_2$ not both zero.  We may rewrite this as
$$\eta_1^{m_1-m_2}(\eta_1\eta_2)^{m_2}=1=\eta_1^{m_1-m_2}t^{-m_2}.$$
Choose $\tau$ such that $\tau^{m_1-m_2}=t$.  Thus $\tau^{m_2}=\eta_1$ and $\tau^{-m_1}=\eta_2$.  It follows that $\tau$ is a root of the following equations:
$$\psi_n(x^{m_1-m_2})=0, {\rm\  and\ \ }x^{m_2}+x^{-m_1}={1\over 1+x^{m_1-m_2}}  \eqno(5.4)$$
We know that $\psi_n(x^{m_1-m_2})$ has exactly $m_1-m_2$ roots outside the closed unit disk, and these are $\lambda_{n}^{1\over m_1-m_2}$,  $\lambda_{n}^{1\over m_1-m_2}\omega$, \dots, $\lambda_{n}^{1\over m_1-m_2}\omega^{m_1-m_2-1}$, where $\omega^{m_1-m_2}=1$. 

Now let $\psi_n(x^{m_1-m_2})=\phi_1(x)\cdots\phi_s(x)$, $\phi_j\in\bZ[x]$ be the irreducible factorization.  We may assume that $\tau$ is a root of $\phi_1$.  Since $\tau$ is not a root of unity, the roots of $\phi_1$ cannot all lie on the unit circle, and thus $\phi_1$ must have a root $\tau_*=\lambda_{n}^{1\over m_1-m_2}\omega^\ell$ for some $\ell$.  Since $\tau$ satisfies the second equation in (5.4), and $\phi_1$ is irreducible, it follows that $\tau_*$ also satisfies this equation, which may be rewritten as
$$\tau_*^{m_2}+{1\over\lambda_n}{1\over \tau_*^{m_2}}={1\over 1+\lambda_n}.$$
Now we write $\tau_*^{m_2}=\delta e^{i\theta}$ with $\delta>0$, $\delta\ne1$, so this equation becomes
$$\delta(\cos\theta+i\sin\theta)+{1\over\lambda_n \delta}(\cos\theta-i\sin\theta)={1\over 1+\lambda_n}.$$
The imaginary part is $\delta\sin\theta-{1\over\lambda_n\delta}\sin\theta=0$, so we must have either $\sin\theta=0$ or $\delta=\lambda_n^{-1/2}$.
The first case is impossible.  In the second case we have $\lambda^{-1/2}_n=\delta=\lambda_n^{m_2/(m_1-m_2)}$, which means that $m_1=-m_2$.   By our condition of multiplicative dependence, we find that $\eta_1=\eta_2$, which implies that  $4t^2+7t+4=0$, which is impossible.

The case $j=2$ is similar.  Suppose again that we have a multiplicative dependence $\eta_1^{m_1}\eta_2^{m_2}=1$, and we choose $\tau$ as before.  Then $\tau$ is a root of
$$\psi_n(x^{m_1-m_2})=0,{\rm\ and \ \ }x^{m_2}+x^{-m_1}={1+x^{m_1-m_2}\over 1+x^{m_1-m_2}+x^{2(m_1-m_2)} }.$$
We now have $\tau_*$ as before, which must satisfy the second equation, which gives
$$\tau_*^{m_2}+{1\over\lambda_n}\tau_*^{-m_2}={1+\lambda_n\over 1+\lambda_n+\lambda_n^2}.$$
Again, we set $\tau_*^{m_2}=\delta e^{i\theta}$, so this equation becomes
$$\delta(\cos\theta+i\sin\theta)+{1\over\lambda_n\delta}(\cos\theta-i\sin\theta) = {1+\lambda_n\over 1+\lambda_n +\lambda_n^2}.$$
Setting the imaginary part equal to zero gives $\sin\theta=0$ or $\delta=\lambda_n^{-1/2}$.  The first possibility is impossible.  The second possibility leads to the conclusion that $\eta_1=\eta_2$, which means that $4t^4+7t^3+10t^2+7t+4=0$.  This is impossible, so $\eta_1$ and $\eta_2$ are multiplicatively independent.
\qed

\proclaim Lemma 5.4.  Suppose that $n\ge8$, and let $j=2$ or 3.  If  $j$ divides $n$, then there exist $(a,b)\in\Gamma_j\cap\cV_n$ such that $|\eta_1|=|\eta_2|=1$.  

\noindent{\it Proof. }  Let us recall the factoriztion $\chi_n=\tilde C_n\psi_n$, with $\tilde C_n$ given as in Theorem 3.5.  Thus the roots of $\chi_n$ are given by $\lambda_n$, $\lambda_n^{-1}$, $1$, the roots of unity $\{r_1,\dots,r_N\}$ coming from $\tilde C_n$, and those roots $\{t_1,\dots,t_M\}$ of $\psi_n$ which have modulus 1.  Since the sum of the roots of $\chi_n$ is zero, we must have
$$\lambda_n+\lambda_n^{-1}+1+Re\sum r_\ell + Re\sum t_\ell =0.$$
On the other hand, $1<\lambda_n<1.4$, and the roots of each cyclotomic factor in $\tilde C_n$ add to either 0 or $-1$.  Thus we have
$$3.4+ Re\sum t_\ell\ge0.$$
Now we claim that we must have a value $t_\ell$ which satisfies both conditions (1) and (2) in Lemma 2.2.  Otherwise, we have $Re( t_\ell)\le -5.8$ for all $\ell$, which means that $M\le 3.4/.58$, so that $M\le 6$.  By Theorem 3.5, the degree of $\psi_n$ is at least $n-26$, so we must have $n\le 34$.  The cases $8\le n\le 34$ may be checked directly.  \qed

 The following result is a combination of Theorems 5.1, 5.3 and Lemma 5.4;  Theorem~10.1 of [M2] gives rank 2 rotation domains under the additional assumption that $n$ is sufficiently large.

 \proclaim Theorem 5.5.  Suppose that $n\ge8$, $j=2$ or 3, and $j$ divides $n$.  Then there exists $(a,b)\in\Gamma_j\cap\cV_n$ such that $f_{a,b}$ has a rank 2 rotation domain centered at $FP_r$, as well as a rank 1 rotation domain centered at $FP_s$.

\noindent{\bf \S6.  Real Mappings of Maximal Entropy. }  Here we consider real parameters $(a,b)\in\bR^2\cap\cV_n$ for $n\ge7$.  Given such $(a,b)$, we let $\cX_R$ denote the closure of $\bR^2$ inside $\cX_{a,b}$. We let $\lambda_n>1$ be the largest root of $\chi_n$, and for $1\le j\le 3$, we let $f_{j,R}$ denote the automorphism of $\cX_R$ obtained by restricting $f_{a,b}$ to $\cX_R$, with $(a,b)=\varphi_j(\lambda_n)$.  
\proclaim Theorem 6.1.  There is a homology class $\eta\in H_1(\cX_R)$ such that $f_{1,R*}\eta=-\lambda_n\eta$.  In particular, $f_{1,R}$ has entropy $\log\lambda_n$. 

\noindent{\it Proof. }  We use an octagon in Figure 6.1 to represent $\cX_R$.  Namely, we start with the real projective plane $\bR{\bf P}^2$; we identify antipodal points in the four ``slanted'' sides.  The  horizontal and vertical pairs of sides of the octagon represent the blowup fibers  over the points $e_1$ and $e_2$.  These are labeled $E_1$ and $E_2$;  the letters along the boundary indicate the identifications.  (Since we are in a blowup fiber, the identification is no longer ``antipodal.'')  Further,  the points $f^jq$ (written ``$j$'') $0\le j\le n$, are blown up, although we do not draw the blowup fibers explicitly.  To see the relative positions of ``$j$'' with respect to the triangle $\Sigma_\beta,\Sigma_\gamma,\Sigma_0$, consult Figure 2.1.  The 1-chains of the homology class $\eta$ are represented by the directed graph $\cG_1$ inside the manifold $\cX_R$.  If we project $\cG_1\subset\cX_R$ down to the projective plane, then all of the incoming arrows at a center of blowup  ``$j$'' will be tangent to each other, as well as the outgoing arrows.  Thus the projection of $\cG_1$ will look like a train track.
\epsfysize=2.8in
\centerline{ \epsfbox{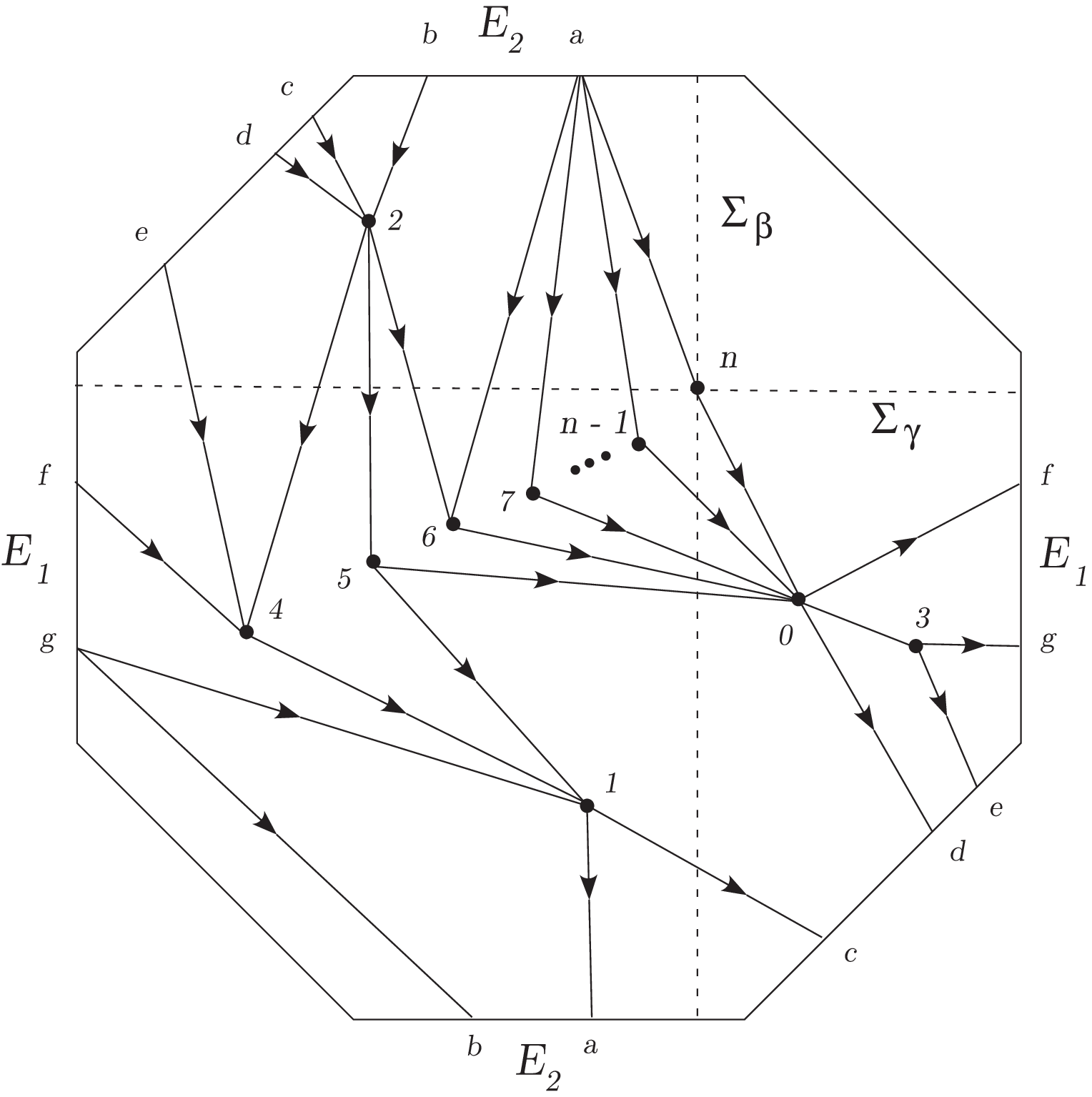}  }

\centerline{Figure 6.1.  Graph $\cG_1$; Invariant homology class for family $\Gamma_1$.}
\smallskip

In order to specify the homology class $\eta$, we need to assign real weights to each edge of the graph.  By ``51'' we denote the edge connecting ``5'' and ``1''; and ``$170=1a70$'' denotes the segment starting at ``1'', passing through ``$a$'', continuing through ``7'', and ending at ``0''.  Abusing notation, also write ``51'', etc., to denote the weight of the edge, as well as the edge itself.  We determine the weights by mapping $\eta$ forward.  We find that, upon mapping by $f$, the orientations of all arcs are reversed.  Let us describe how to do this.  Consider the arc ``34''=``$3e4$''.  The point ``$e$'' belongs to $\Sigma_0$, and so it maps to $E_1$.  Thus ``34'' is mapped to something starting at ``4'', passing through $E_1$, and then continuing to ``5''.  Thus we see that ``34'' is mapped (up to homotopy) to ``4f05''.  Thus, the image of ``34'' covers ``04'' and ``05''.

Inspection shows that no other arc maps across ``04'', so we write ``$04\to  34$'' to indicate that the weight of side ``04'' in $f_*\eta$ is equal to the weight of ``34''.  Inspecting the images of all the arcs, we find that ``24'' also maps across ``05'', so we write ``$05\to  24 + 34$'' to indicate how the weights transform as we push $\cG_1$ forward.   Looking at all possible arcs, we write the transformation $\eta\mapsto f_*\eta$ as follows:
$$\eqalign{&02\to  16 + 170 + \cdots + 1(n-1)0, \ \ \ \ 03\to  24 + 25 +26, \ \ \ \ 04\to  34,\cr
& 05\to  24 + 34,\ \ \ \ 
06\to  25, \ \ \ \ 12\to  1n0, \ \ \ \ 13\to  02,\ \ \ \ 14\to  03,\cr
& 15\to  04,\ \ \ \   16\to  05,\ \ \ \ 170\to  06,\ \ \ \ 
1k0\to  1(k-1)0,\  7<k\le n,\cr
&23\to  12, \ \ \ \ 24\to  13, \ \ \ \   25\to  14, \ \ \ \ 26\to  15, \ \ \ \ 34\to  23,\cr}\eqno(6.1)$$

The formula (6.1) defines a linear transformation on the space of coefficients of the 1-chains defining $\eta$.   The spectral radius of the transformation (6.1) is computed in Appendix C, where we find that it is $\lambda_n$.  Now let $w$ denote the eigenvector of weights corresponding to the eigenvalue $\lambda_n$.  It follows that if we assign these weights to $\eta$, then by construction we have $f_{1,R*}\eta=-\lambda_n\eta$, and $\eta$ is closed.  \qed 

\noindent{\bf Remark. }  Let us compare with the situation for real H\'enon maps.  In [BLS] it was shown that a real H\'enon map has maximal entropy if and only if all periodic points are real.  On the other hand, if $(a,b)=\varphi_1(t)$, $1\le t\le 2$,  the (unique) 2-cycle of the map $f_{a,b}$ is non-real.  This includes all the maps discussed in Theorem 6.1, since all values of $t=\lambda_n$ are in this interval.

\epsfysize=2.8in
\centerline{ \epsfbox{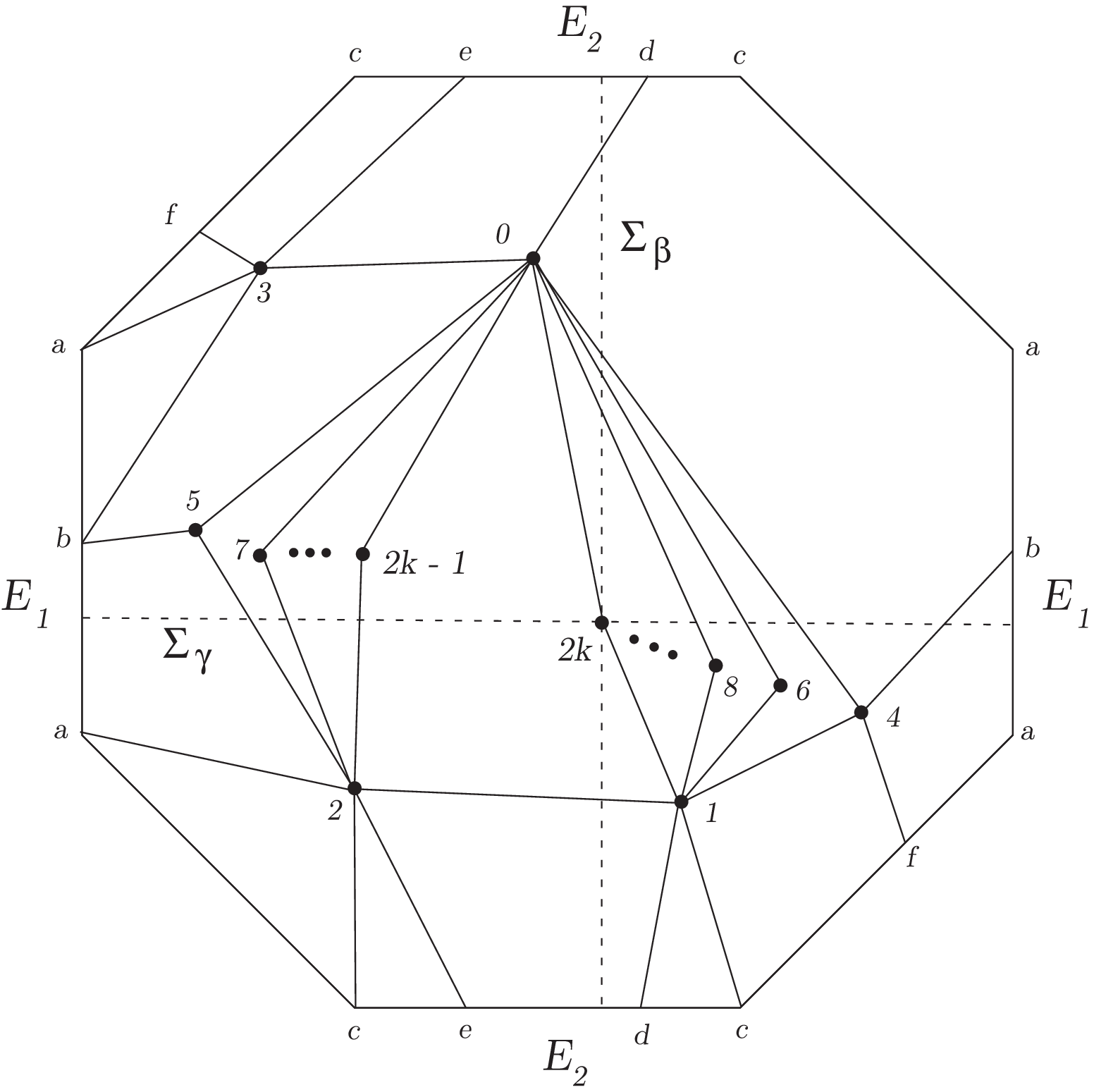}  }

\centerline{Figure 6.2.  Invariant graph $\cG_2$: $n=2k$.}

\proclaim Theorem 6.2. The maps $f_{2,R}$ (if $n$ is even) and  $f_{3,R}$ (if $n$ is divisible by 3) have entropy $\log\lambda_n$.

\noindent{\it Proof. }    Since the entropy of the complex map $f$ on $\cX$ is $\log\lambda_n$, the entropies of $f_{2,R}$ and $f_{3,R}$ are bounded above by $\log\lambda_n$.  In order to show that equality holds for the entropy of the real maps, it suffices by Yomdin's Theorem [Y] (see also [G]) to show that $f_{j,R}$ expands lengths by an asymptotic factor of $\lambda_n$.  We will do this by producing graphs $\cG_2$ and $\cG_3$ on which $f$ has this expansion factor.
We start with the case $n=2k$; the graph $\cG_2$ is shown in Figure 6.2, which should be compared with Figure 2.2.   We use the notation $01=0d1$ for the edge in $\cG_2$ connecting ``0'' to ``1'' by passing through $d$.  In this case, the notation already defines the edge uniquely; we have added the $d$ by way of explanation.  We will use $\cG_2$ to measure length growth; the edges of $\cG_2$ will be analogous to the ``u-arcs'' that were used in [BD2].   By ``01'' we will mean all arcs in $\cX_R$ which connect the fiber ``0'' to the fiber ``1``, and whose projections to $\bR\bP^2$ are homotopic to the edge 01.  In order to measure length growth, we will start with a configuration of arcs corresponding to edges of the graph.  The symbol ``01'' will also be used in formulas (6.2) and (6.3), for instance, to denote the number of arcs which are of type 01.

Now we discuss how these arcs are mapped.  The arc $01$ crosses $\Sigma_\beta$ and then $E_2\ni d$ before continuing to ``1''.  Since $\Sigma_\beta$ is mapped to $E_2$ and $E_2$ is mapped to $\Sigma_0$, the image of 01 will start at ``1'' and cross $E_2$ and then $\Sigma_0$ before reaching ``2''.   Up to homotopy, we may slide the intersection points in  $E_2$ and $\Sigma_0$ over to a point $g\in E_2\cap\Sigma_0$.  Thus, up to homotopy, $f$ maps the edge $01$ in $\cG_2$  to the edge $1g2$.
\smallskip
\epsfysize=2.8in
\centerline{ \epsfbox{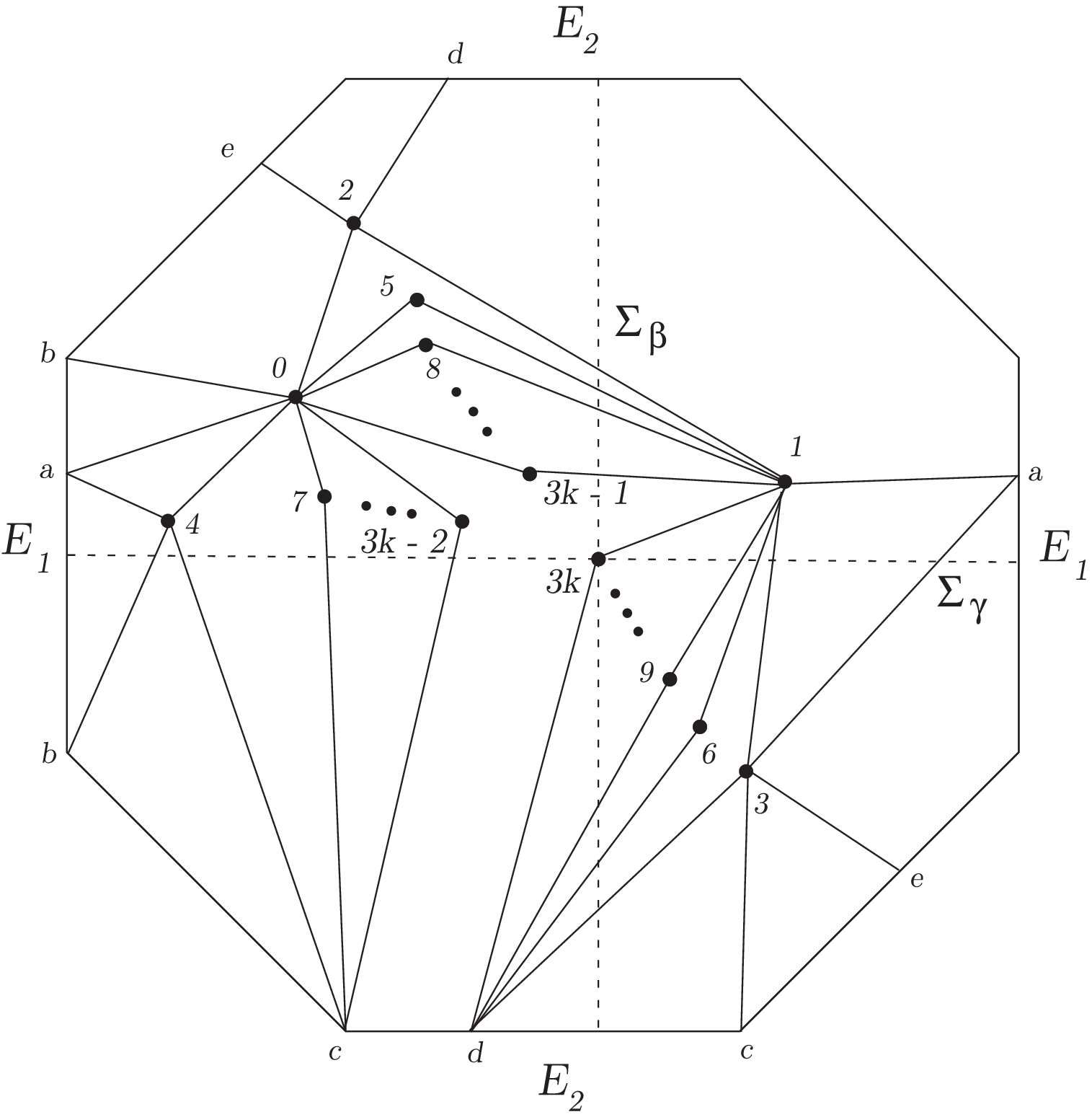}  }

\centerline{Figure 6.3.  Invariant graph $\cG_3$: $n=3k$.}
\smallskip

Similarly, we see that the arcs $04$, $06$, \dots, $0(2k-1)$ all cross $\Sigma_\beta$ and then $\Sigma_\gamma$.  Thus the images of all these arcs will start at $1$, pass through $E_2$ at $d$, then $0$, and continue to the respective endpoints $5$, $7$,  \dots, $(2k-1)$.  Since the images of all these arcs, up to homotopy, contain the edge $01$, the transformation of weights in the graph is given by the first entry of (6.2), and the whole transformation is given by the rest of (6.2):
$$\eqalign{
&01\to  04+061 + 081 +\cdots+0(2k-1)1,\ \ \ \ 03\to  25 + 072 + 092+\cdots + 0(2k-1)2\cr
&04\to  3b4,\ \ \ \ 061\to  25 + 45, \ \ \ \  12\to  0(2k)1, \ \ \ \ 1c2\to  01,\ \ \  14\to  03\cr 
& 0(2j)1\to  0(2j-1)2,\ j= 4,5,\dots,k, \ \ \ \ 0(2j+1)2\to  0(2j)1, \ \  j=3,4,\dots,k-1\cr
&2a3\to  12, \ \ \ \ 2e3\to  1c2,\ \ \ \ 25\to 14, \ \ \ \ 3b4\to  2e3, \ \ \ \ 34\to  2a3, \ \ \ \ 45\to  34.\cr}\eqno(6.2)$$
The characteristic polynomial for the transformation defined in (6.2) is computed in the Appendix, and the largest eigenvalue of (6.2) is $\lambda_n$, so $f_{2,R}$ has the desired expansion.

The case $n=3k$ is similar.  The graph $\cG_3$ is given in Figure 6.3.  Up to homotopy, $f_{3,R}$ maps the graph $\cG_3$ to itself according to:
$$\eqalign{
&01\to  0b4, \ \ \ \ 02\to  13+162+192+\cdots+1(3k-3)2,\cr
&0b4\to  3a4+073+0(10)3+\cdots +0(3k-2)3,\ \ \ \ 051\to  0b4+3c4+3a4,\cr
&0(3j-1)1\to  0(3j-2)3, \ j=3,4,\dots,k,
\ \ \ \  1(3j)2\to  0(3j-1)1,\ j=2,3,\dots,k,\cr
&0(3j+1)3\to  1(3j)2, \ j=2,3,\dots,k-1,\cr
&12\to 01,\ \ \ \ 13\to 02, \ \ \ \ 23\to  1(3k)2, \ \ \ \ 2d3\to 12, \ \ \ \ 3a4\to  23, \ \ \ \ 3c4\to  2d3.
}\eqno(6.3)$$
The linear transformation corresponding to (6.3) is shown in Appendix C to have spectral radius equal to  $\lambda_n$, so $f_{3,R}$ has entropy $\log\lambda_n$.  \qed

\bigskip\noindent{\bf\S7.  Proof of  Theorem B. }  This Theorem is a consequence of results we have proved already.  Let $f_{a,b}$ be of the form (0.1).  By Proposition B.1 in Appendix B, we may suppose that $n\ge11$.  Thus if $f$ has an invariant curve, then by Theorem 1.2, it has an invariant cubic, which is given explicitly by Theorem 2.1.  Further, by Theorem 2.2, we must have $(a,b)=\varphi_j(t)$ for some $j$ dividing $n$, and a value $t\in\bC$ which is a root of $\chi_n$.  By Theorem 3.3, $t$ cannot be a root of unity.  Thus it is a Galois conjugate of $\lambda_n$.  The Galois conjugates of $\lambda_n$ are of two forms: either $t$ is equal to $\lambda_n$ or $\lambda_n^{-1}$, or $t$ has modulus equal to 1.  

In the first case, $(a,b)\in\cV\cap\bR^2$, and thus $f$ is a real mapping.  The three possibilities are $(a,b)\in\cV_n\cap\Gamma_j$, $j=1,2,3$, and these are treated in \S5.  In all cases, we find that the entropy of the real mapping $f_{R,j}$ has entropy equal to $\log\lambda_n$.  By Cantat [C2], there is a unique measure $\mu$ of maximal entropy for the complex mapping.  Since $f_{R,j}$ has a measure $\nu$ of entropy $\log\lambda_n$, it follows that $\mu=\nu$, and thus $\mu$ is supported on the real points.  On the other hand, we know that $\mu$ is disjoint from the Fatou sets of $f$ and $f^{-1}$.  McMullen [M2] has shown that the complement of one of the Fatou sets $\cF(f)$ or $\cF(f^{-1})$ has zero volume.  The same argument shows that the complement  inside $\bR^2$ has zero area.  Thus the support of $\mu$ has zero planar area.

The other possibility is that $t$ has modulus 1.  In this case, Theorem B is a consequence Theorem 5.1.  \qed

\bigskip\noindent{\bf \S8.  Representing the Coxeter Element. }  Now we suppose that ${\cal X}$ is a rational surface, and $F\in Aut({\cal X})$ is an automorphism such that the induced map $F_*$ on $Pic({\cal X})$ has infinite order.  Then by a Theorem of Nagata there is a biholomorphic morphism $\pi:{\cal X}\to{\bf P}^2$ such that $\pi$ can be factored $\pi_1\circ\cdots\circ\pi_N$ into a composition of point blow-ups $\pi_\sigma:{\cal X}_\sigma\to{\cal X}_{\sigma-1}$ such that ${\cal X}_N={\cal X}$ and ${\cal X}_0={\bf P}^2$.  For each $1\le\sigma\le N$ we let $E_\sigma$ denote the exceptional blowup fiber for $\pi_\sigma$.  If we let $E_0$ denote the class of a general line, then $\{E_\sigma:0\le\sigma\le N\}$ is a basis of $Pic({\cal X})$.  We may also define ${\cal E}_\sigma:=(\pi_{\sigma+1}\circ\cdots\circ\pi_N)^*E_\sigma$, and ${\cal E}_0=E_0$.  Thus $\{{\cal E}_\sigma:0\le\sigma\le N\}$ is a geometric basis of $Pic({\cal X})$, which means that ${\cal E}_0^2=1$, ${\cal E}^2_\sigma=-1$ for $1\le\sigma\le N$, and ${\cal E}_i\cdot {\cal E}_j=0$ when $i\ne j$.

Another result of Nagata says that $F_*$ must belong to the Weyl group $W_N$, which is generated by reflections $r_0,\dots,r_{N-1}$.  If $1\le\sigma\le N-1$, $r_\sigma:={\cal E}_\sigma\leftrightarrow {\cal E}_{\sigma+1}$ is the reflection about the element ${\cal E}_{\sigma}-{\cal E}_{\sigma+1}$; and $r_0$ is the Cremona involution, which is the reflection about ${\cal E}_0-{\cal E}_1-{\cal E}_2-{\cal E}_3$.   The Coxeter is the product $r_1\cdots r_{N-1} r_0$ of these involutions.  Writing the product $r_1\cdots r_{N-1}$ in a different order corresponds to permuting the basis elements ${\cal E}_\sigma$, $1\le \sigma\le N$.  Thus we may write the Coxeter element (after a possible permutation of basis elements) as:
$$\eqalign{ {\cal E}_0 &\mapsto  2{\cal E}_0 - {\cal E}_2 - {\cal E}_3 - {\cal E}_4 \cr
  {\cal E}_1 &\mapsto  {\cal E}_0 -{\cal E}_3 - {\cal E}_4\cr
  {\cal E}_2 & \mapsto  {\cal E}_0 - {\cal E}_2 - {\cal E}_4\cr
  {\cal E}_3 & \mapsto {\cal E}_0 - {\cal E}_2 - {\cal E}_3\cr
  {\cal E}_\sigma & \mapsto {\cal E}_{\sigma +1}, \ \ 4\le\sigma\le N-1, \ \ \ {\cal E}_N\mapsto{\cal E}_1.}$$
\proclaim Theorem 8.1.  Suppose that ${\cal X}$ is a rational surface,  $F\in Aut({\cal X})$ represents the Coxeter element of $W_N$ for $N\ge 5$, and suppose that we represent ${\cal X}$ by iterated blowups $\pi:{\cal X}\to{\bf P}^2$.  Then there exist a linear map $T\in Aut({\bf P}^2)$ and complex numbers $a$ and $b$ such that
$$f_{a,b}\circ T\circ\pi = T\circ \pi\circ F.$$

Before giving the proof, let us recall that for $\sigma\ge1$, the ${\cal E}_\sigma$ are given as $E_\sigma+\sum_t \mu_{\sigma, t} E_t$, where $\mu_{\sigma,t}$ is a  nonnegative integer, and $E_t$ is an irreducible exceptional fiber ``over''  $E_\sigma$.  In addition, each of the basis elements ${\cal E}_\sigma$, $1\le \sigma\le N$ lies over a single point $p_\sigma:=\pi(E_\sigma)$.  It is uniquely represented in terms of prime components, so the number $\#{\cal E}_\sigma$ of these prime components is well defined for $1\le\sigma\le N$.

\proclaim Lemma 8.2.  We have $\#{\cal E}_\sigma=1$, i.e., ${\cal E}_\sigma=E_\sigma$ for $4\le\sigma\le N$, and $\sigma=1$.

\noindent{\it Proof. }  Since the representation $E_\sigma+\sum_t \mu_{\sigma, t} E_t$ is unique, we have $\#{\cal E}_4=\cdots=\#{\cal E}_N=\#{\cal E}_1$.  Now if $\#{\cal E}_4>1$, then there is some $\tau=\tau_4$ such that ${\cal E}_{\tau_4}$ lies strictly above ${\cal E}_4$, and thus $\#{\cal E}_4>\#{\cal E}_{\tau_4}\ge1$.  Now $F_*$ maps ${\cal E}_{\tau_4}$ to some basis element ${\cal E}_{\tau_5}$, which lies strictly above ${\cal E}_{5}$.  It is evident that we have $\#{\cal E}_{\tau_4}=\#{\cal E}_{\tau_5}=\#{\cal E}_{\tau_1}<\#{\cal E}_4$, so we have three new independent elements.  But now we have too many independent elements for the dimension of $Pic({\cal X})$, so we conclude that we must have had $\#{\cal E}_4=1$. \qed

\proclaim Lemma 8.3.  The points $p_2$ and $p_3$ are distinct.

\noindent{\it Proof. }  The Coxeter element above maps ${\cal E}_3\mapsto {\cal E}_0-{\cal E}_2-{\cal E}_3\mapsto {\cal E}_2$, so $\#{\cal E}_2=\#{\cal E}_3$.  If $\#{\cal E}_3=\#{\cal E}_2=1$, then we have $\#{\cal E}_\sigma=1$ for all $1\le \sigma\le N$. Thus the points $p_\sigma$ are distinct.  Now suppose that $\#{\cal E}_3=\#{\cal E}_2>1$.  Then in the representation ${\cal E}_3=E_3+\sum_t \mu_{3, t} E_t$ there are components $E_t$ lying above $E_3$.  However, if $p_2=p_3$, then either $E_2$ lies above $E_3$ or $E_3$ lies above $E_2$.  In the first case, for instance, we see from this representation that we must have  $\#{\cal E}_2<\#{\cal E}_3$.  Thus we must have  $p_2\ne p_3$.
\qed

\proclaim Lemma 8.4.  The points $p_1$, $p_2$, $p_3$ and $p_4$ are distinct

\noindent{\it Proof. }  Let us start by showing that various of the $p_i$'s are distinct.  First we show $p_1\ne p_2$.  Suppose not.  If $p_1=p_2$, then since $\#{\cal E}_1=1$, we have ${\cal E}_1=E_1$, and the fiber $E_1$ must lie above $E_2$.  And in the representation ${\cal E}_2=E_2+E_1 +\sum_{t\ne 1,2,3}\mu_{2,t}E_t$, and all the other terms $\mu_{2,t}$ are either 0 or 1.   Thus
$$E_2={\cal E}_2-{\cal E}_1 -\sum_{t\ne 1,2,3}\mu_{2,t}E_t,$$
and mapping this expression forward by the Coxeter element gives
$$E_2\mapsto   {\cal E}_3-{\cal E}_2 - \sum_{t\ne 1,2,3}\mu_{2,t}E_{t+1}.$$
But this is not possible, since $E_2$ is positive and irreducible, but the image in the line above is neither positive nor irreducible, since ${\cal E}_2$ and ${\cal E}_3$ are carried over distinct points $p_2\ne p_3$.

Similarly, we show that $p_2\ne p_4$.  For if $p_2=p_4$, then we must have $E_4$ above $E_2$, and in the expression ${\cal E}_2 = E_2+\sum_t\mu_{2,t}E_t$, we have $\mu_{2,4}=1$.  The inverse of the Coxeter element acts like ${\cal E}_2\mapsto {\cal E}_0-{\cal E}_2-{\cal E}_3$ and ${\cal E}_4\mapsto{\cal E}_0-{\cal E}_1-{\cal E}_2$.  Thus the inverse acts by
$$E_2\mapsto  {\cal E}_1-{\cal E}_3 - \sum_{t\ne 2,3,4}\mu_{2,t}E_{t-1},$$
which is not possible since the right hand side is not positive since ${\cal E}_3-{\cal E}_1\ge E_3$.

Next we show that $p_1\ne p_3$ and $p_4\ne p_3$ by the same argument.  Finally, we see that $p_1\ne p_4$.  Since $p_1$ and $p_4$ are disjoint from $\{p_2,p_3\}$, it follows that there is just one level of blowup in ${\cal X}$ over $p_1$ and $p_4$.  Since $E_\sigma={\cal E}_\sigma$ for $\sigma=1,4$, and  these are distinct elements of $Pic({\cal X})$, $p_1\ne p_4$.
\qed

\proclaim Lemma 8.5.  $p_1,p_4\notin \overline{p_2p_3}$.

\noindent{\it Proof. }  By Lemma 8.2 the complete fiber over $p_1$ is $E_1={\cal E}_1$, which is irreducible.  Thus the image $F(E_1)$ is an irreducible curve in ${\cal X}$ which represents ${\cal E}_0-{\cal E}_3-{\cal E}_4$.  The curve $\pi(F(E_1))$ is a curve in ${\bf P}^2$ of degree 1 and thus is a line which must pass through $p_3$ and $p_4$.   Since $F(E_1)$ is irreducible, it cannot pass through $p_1$ or $p_2$.  In particular, the points $p_2,p_3,p_4$ are not colinear, so $p_4\notin \overline{p_2p_3}$.  Similarly, we apply the inverse of the Coxeter element fo $E_4={\cal E}_4$ and conclude that $p_3$ and $p_4$ are not on the line through $p_1p_2$.   Thus $p_1$, $p_2$ and $p_3$ are not collinear, so $p_1\notin \overline{p_2p_3}$.  \qed

\noindent{\it Proof of Theorem 8.1. }  After Lemmas 8.2--5, we are in the situation of Theorem 11.1 of [M2].  We give the proof here for the sake of completeness.  Let $\varphi$ be the birational map of ${\bf P}^2$ induced by $F$.  Since $p_1\notin\overline{p_2p_3}$, the points $p_1$, $p_2$ and $p_3$ are the vertices of a proper triangle.  Thus we may introduce a linear change of coordinates $[t:x:y]$ such that $p_1=[1:0:0]$, $p_2=[0:1:0]$, and $p_3=[0:0:1]$.  Now define $\alpha$ and $\beta$ by the condition that $p_4=[1:\alpha:\beta]$.  Define $g_{\alpha,\beta}= (y,y/x) + (\alpha,\beta)$.  Then since $\varphi$ and $g_{\alpha,\beta}$ have the same exceptional curves and map them in the same way,  there are constants $A$ and $B$ such that $\varphi(Ax,By)=g_{\alpha,\beta}(x,y)$.  Since $g_{\alpha,\beta}$ is linearly conjugate to some $f_{a,b}$, the Theorem follows.
\qed

Let $\pi:{\cal Y}\to{\bf P}^2$ denote the space obtained by blowing up the points $e_1$ and $e_2$ in ${\bf P}^2$.  Every map $f_{a,b}$ induces a birational map $f_{\cal Y}$ of ${\cal Y}$ to itself, and $p=(-b,-a)$ is the only point of indeterminacy for $f_{\cal Y}$.  Thus the pointwise iteration $f^j_{\cal Y}:=f_{\cal Y}\circ\cdots\circ f_{\cal Y}$ is well defined until the orbit reaches $p$.  In other words:
\proclaim Lemma 8.6.  If $r\in{\cal Y}$ and $f^j_{\cal Y}r\ne p$ for $0\le j\le k$, then $(f^k)_{\cal Y}(r)=f_{\cal Y}\circ\cdots\circ f_{\cal Y}(r)$, there the composition on the right is $k$-fold.

The point of Lemma 8.6 is that since we may identify ${\bf C}^2$ with ${\cal Y}-(E_1\cup E_2\cup\Sigma_0)$, we may define ${\cal V}_n$ in terms of $f_{\cal Y}$:
$${\cal V}_n=\{(a,b)\in {\bf C}^2: f_{\cal Y}^jq\ne p{\rm\ for\ }0\le j<n, \ f_{\cal Y}^nq=p\}.$$
By Lemma 8.2, a representation of the Coxeter element can be obtained by making simple (not iterated) blowups of birational maps $f_{a,b}$ which have been lifted to ${\cal Y}$.  We state it as follows:
\proclaim Theorem 8.7.  Suppose that ${\cal X}$ is a rational surface, and $F\in Aut({\cal X})$ represents the Coxeter element in $W_N$.  Then there exists a surface $\pi:{\cal Z}\to {\cal Y}$ which is obtained by the blowup of a finite set of distinct points in ${\cal Y}$, and there are an automorphism $G\in Aut({\cal Z})$ and $(a,b)\in {\cal V}_{N-3}$ such that $(F,{\cal X})$ is conjugate to $(G,{\cal Z})$, and $\pi\circ G=f_{\cal Y}\circ \pi$, with $f_{\cal Y}=(f_{a,b})_{\cal Y}$.

\noindent {\bf Remarks. } Theorems 8.1 and 8.7  give a strengthening of Theorem 11.1 of [M2] so that it applies to all rational surface maps.   We note some notational differences with [M2].  Namely, [M2] uses $N$ for the number of blowups, whereas [BK] and the present paper use $n$ for the length of the orbit of $q$.  These two numbers are related by $N=n+3$.  The 2-dimensional surface on which an automorphism is defined is denoted by ${\cal X}_{a,b}$ in [BK], whereas the surface is $S$ in [M2], and $X$ represents an invariant curve.  Here we use the notation (0.1) of [BK] for the birational family; [M2] changed this notation to $(x,y)\mapsto (y,y/x)+(\alpha,\beta)$, which is conjugate to $f_{a,b}$ if $\alpha=b-a$ and $\beta=a$.

\bigskip\noindent{\bf Appendix A.  Varieties $\cV_j$ and $\Gamma$ for $0\le j\le 6$. }  The sets $\cV_j$, $0\le j\le 6$ are enumerated in [BK].  We note that $\cV_0=(0,0)\subset\Gamma_2\cap\Gamma_3$, $\cV_1=(1,0)\subset \Gamma_1$, $\cV_2\subset\Gamma_1\cap\Gamma_2$, $\cV_3\subset\Gamma_1\cap\Gamma_3$, and $(\cV_4\cup\cV_5)\cap\Gamma=\emptyset$.    

Each of the mappings in $\cV_6$ has an invariant pencil of cubics.  There are two cases: the set $\cV_6\cap\{b\ne0\}$ (consisting of four points) is contained in $\Gamma_2\cap\Gamma_3$.
The other case, $\cV_6\cap\{b=0\}=\{(a,0): a\ne0,1\}$, differs from the cases $\cV_n$, $n\ne 0,1,6$, because the manifold $\cX_{a,0}$ is constructed by iterated blowups ($f^4q\in E_1$ and $f^2q\in E_2$, see [BK, Figure~ 6.2]).  

The invariant function $r(x,y)=(x+y+a)(x+1)(y+1)/(xy)$ for $f_{a,0}$, which defines the invariant pencil, was found by Lyness [L]  (see also [KLR],  [KL], [BC] and [Z1]).  We briefly describe the behavior of $f_{a,0}$. 
By $M_\kappa=\{r=\kappa\}$ we denote the level set of $r$ inside $\cX_{a,0}$.  The curve $M_\infty$ consists of an invariant 5-cycle of curves with self-intersection $-2$:
$$\Sigma_\beta=\{x=0\}\mapsto E_2\mapsto\Sigma_0\mapsto E_1\mapsto \Sigma_B=\{y=0\}\mapsto\Sigma_\beta.$$
The restriction of $f^5$ to any of these curves is a linear (fractional) transformation, with multipliers $\{a,a^{-1}\}$ at the fixed points.   $M_0$ consists of a 3-cycle of curves with self-intersection $-1$:
$\{y+1=0\}\mapsto\{x+1=0\}\mapsto\{x+y+a=0\}.$  The restriction of $f^3$ to any of these lines is linear (fractional) with multipliers $\{a-1,(a-1)^{-1}\}$ at the fixed points.

\proclaim Theorem A.1.  Suppose that $a\notin\{-{1\over 4},0,{3\over 4}, 1,2\}$, and $\kappa\ne0,\infty$.  If $M_\kappa$ contains no fixed point, then $M_\kappa$ is a nonsingular elliptic curve, and $f$ acts as translation on $M_\kappa$.  If $M_{\kappa}$ contains a fixed point $p$, then $M_{\kappa}$ has a node at $p$.  If we uniformize $s:\hat\bC\to M_{\kappa}$ so that $s(0)=s(\infty)=p$, then $f|_{M_{\kappa}}$ is conjugate to $\zeta\mapsto\alpha\zeta$ for some $\alpha\in\bC^*$.

The intersection $\Gamma_j\cap\{b=0\}\cap\cV_6$ is given by $(-{1\over 4},0)$, $({3\over4},0)$, or $(0,2)$, if $j=1,2$, or 3, respectively.

\proclaim Theorem A.2.  Suppose that $a=-{1\over 4},{3\over4},$ or 2.  Then the conclusions of Theorem A.1 hold, with the following exception.  If $FP_s\in M_{\kappa}$, then $M_\kappa$ is a cubic which has  a cusp at $FP_s$, or is a line and a quadratic tangent at $FP_s$, or consists of three lines passing through $FP_s$.  If we uniformize a component $s:\hat\bC\to M_\kappa$ such that $s(\infty)=FP_s$, then $f^j|_{M_\kappa}$ is conjugate to $\zeta\mapsto\zeta+1$, where $j$ is chosen so that $(a,0)\in\Gamma_j$.

\bigskip\noindent{\bf Appendix B.  Varieties $\cV_n$ and $\Gamma$ for $n\ge7$. }  We may compute  $\cV_n$ explicitly by starting with the equation $f_{a,b}^n(-a,0)+(b,a)=0$.  We look upon $f^n_{a,b}(-a,0)$ as a rational expression in $a$ and $b$ with integer coefficients, so we first remove common factors from the numerator and denominator of $f^n$, and then we convert this equation to a pair of polynomial equations $S_n=\{P_n(a,b)=Q_n(a,b)=0\}$ with integer coefficients in the variables $a$ and $b$.  Thus $S_n$ contains all the elements of ${\cal V}_n$ and possibly more.  We may find the elements of $S_n$ by applying elimination theory to $P_n$ and $Q_n$; since they are polynomials with integer coefficients, we obtain an exact polynomial for the possible values of $a$ (or $b$).  We may then test numerically whether a given element $(a,b)\in S_n$ actually belongs to ${\cal V}_n$, and not some ${\cal V}_j$ for $j<n$, by computing the orbit $f_{a,b}^jq$, $1\le j\le n$.  We do this using {\sl Mathematica} (or {\sl Maple}) and find that $\#\cV_7= 10$.  On the other hand, by Theorem 3.5, $\chi_7(x)=(x-1)\psi_7(x)$, and so $\psi_7$ has degree $10$.  By Theorem 3.4, $\#(\cV_7\cap\Gamma_1)=10$, and so by counting we conclude that $\cV_7\subset\Gamma_1$.  Arguing in this manner, we obtain
\proclaim  Proposition B.1.  $\cV_7\subset\Gamma_1$, $\cV_8\cup\cV_{10}\subset\Gamma_1\cup\Gamma_2$, and $\cV_9\subset\Gamma_1\cup\Gamma_3$.

\noindent{\bf Counterexamples. } To find the counterexamples discussed in \S4, we let $n=11$.  By Theorem 3.5, $\psi_{11}$ has degree 10, so by Theorem 3.3, $\Gamma_1\cap\cV_{11}=\Gamma\cap\cV_{11}$ contains 10 elements.  On the other hand, we may use {\sl Mathematica} to determine the number of elements of $\cV_{11}$.  That is, we follow the scheme described in the previous paragraph and obtain the resultant of $P_{11}$ and $Q_{11}$, eliminating either $a$ or $b$.  This resultant has a factor of degree 10 and one of degree 12.  The 10 roots of the degree 10 factor correspond to the values of $\Gamma\cap{\cV}_{11}$, and the 12 roots of the degree 12 factor correspond to the values of ${\cV}_{11}-\Gamma$, and thus these 12 maps do not have invariant curves.

We have repeated this computational procedure in the cases $11<n\le 25$, and we have found that $\cV_n\cap\Gamma$ becomes a quite small fraction of $\cV_n$ as $n$ increases.

\medskip
A map $f_{a,b}$ has a unique 2-cycle and a unique 3-cycle.  For $\ell=2,3$, we let $J_\ell$ denote the product of the Jacobian matrix around the $\ell$-cycle.
\proclaim Theorem B.2.  For $\ell=2,3$, the determinant $\mu_\ell$ of $J_\ell$ is given by
$$\mu_2={a-b-1\over 2b^2+a-1},\ \ \mu_3={1+b+b^2-a-ab\over 1-a-ab}\eqno(B.1)$$
and the trace $\tau_\ell$ is given by
$$\tau_2={3-2a+b-b^2\over 2b^2+a-1},\ \ \tau_3={2+a^2+b+2b^2-b^3+b^4+a(-2-b+2b^2) \over -1+a-ab}.  \eqno(B.2)$$

\noindent{\it Proof. }  The proof of the 2-cycle case is simpler and omitted.  We may identify a 3-cycle with a triple of numbers $z_1,z_2,z_3$:
$$\zeta_1=(z_1,z_2)\mapsto \zeta_2=(z_2,z_3={a+z_2\over b+z_1}) \mapsto \zeta_3=(z_3,z_1={a+z_3\over b+z_2})\mapsto\zeta_1=(z_1,z_2={a+z_1\over b+z_3}).$$
Substituting into this 3-cycle, we find that $z_1,z_2,z_3$ are the three roots of 
$$P_3(z)=z^3+(1+a+b+b^2)z^2+(b^3+ab+2a-1)z-1+a-b+ab-b^2.$$
It follows that 
$$\eqalign{z_1+z_2+z_3&=-(1+a+b+b^2)\cr
z_1z_2+z_1z_3+z_2z_3 &=-1+2a+ab+b^3\cr
z_1z_2z_3 &=1-a+b-ab+b^2\cr}\eqno(B.3)$$
Since
$$Df_{a,b}(\zeta_1=(z_1,z_2))=\pmatrix{0&1\cr -{a+z_2\over (b+z_1)^2} & {1\over b+z_1}\cr} = \pmatrix{0&1\cr {-z_3\over b+z_1}& {1\over b+z_1}\cr},$$
the determinant of $Df_{a,b}(\zeta_1)=z_3/(b+z_1)$, and therefore
$$\mu_3={z_2\over b+z_3}{z_1\over b+z_2}{z_3\over b+z_1}={z_1z_2z_3\over (b+z_3)(b+z_2)(b+z_1)}.$$
Using equations (B.3) we see that $(b+z_1)(b+z_2)(b+z_3)=1-a+ab$ so $\mu_3$ has the form given in (B.1).

Similarly, we compute
$$Tr(J_3)=-{-1+b(z_1+z_2+z_3) +z_1^2+z_2^2+z_3^2\over (b+z_1)(b+z_2)(b+z_3)}.$$
Using (B.3) again, we find that $\tau_3$ is given in the form (B.2). \qed

A computation shows the following:
\proclaim Theorem B.3.  Suppose that $(a,b)=\varphi_1(t)\in\cV_n\cap\Gamma_1$, $n\ge7$.  Then the eigenvalues of $Df$ at $FP_r$ are given by $\{\eta_1=1/t,\eta_2=-(t^3+t^2-1)/(t^4-t^2-t)\}$, where $t$ is a root of $\psi_n$.  Further, they have the resonance $\eta_1^n\eta_2=1$.

In the previous theorem, the fixed point $FP_r$ is contained in the invariant curve.  If $\ell$ divides $n$, for $\ell=2$ or 3, then the $\ell$-cycle is disjoint from the invariant curve.  Thus we have:
\proclaim Theorem B.4.  Suppose $\ell=2$ or $3$ and $n=k\ell \ge7$.  If $(a,b)=\varphi_\ell(t)\in\cV_n\cap\Gamma_\ell$, then the eigenvalues of the $\ell$-cycle are $\{\eta_1=t^{-\ell},\eta_2=-t^{\ell-1}(t^3+t^2-1)/(t^3-t-1)\}$.  Further, they have the resonance $\eta_1^{n+1}\eta_2=1$.

\proclaim Corollary B.5.  If $t=\lambda_n$ or $\lambda_n^{-1}$, then the cycles discussed in Theorems B.3 and B.4 are saddles.  If $t$ has modulus 1, then the multipliers over these cycles have modulus 1 but are not roots of unity. 
\vfill\eject

\bigskip\noindent{\bf Appendix C.  Computation of  Characteristic Polynomials. } 

\proclaim Theorem C.1.  If $\chi_n$ is as in (0.3), and $n\ge7$,  then
\item{(i)}  The characteristic polynomial for  (6.1) is  $(x^7+1)\chi_n(x)/(x^2-1)$; 
\item{(ii)} The characteristic polynomial for (6.2) is $(x^5-1)\chi_{2k}(x)/(x^2-1)$; 
\item{(iii)} The characteristic polynomial for (6.3) is $(x^4-1)\chi_{3k}(x)/(x^3-1)$.

\noindent{\it Proof. }  We start with case $(i)$.  Since the case $n=7$ is easily checked directly, it suffices to prove $(i)$ for $n\ge8$.  Let us use the ordered basis: 
$$\{12, 23, 34, 04, 15, 26, 03, 14, 25, 05, 16, 02, 13, 24, 06, 170, 180, \dots, 1n0\},$$
and let $M=(m_{i,j})$ denote the matrix which represents the transformation $\eta\mapsto f_*\eta$ defined in (5.1), i.e., we set $m_{i,j}=1$ if the $i$-th basis element in our ordered basis maps to the $j$-th basis element, and 0 otherwise.  To compute the characteristic polynomial of $M$, we expand $\det(M-xI)$  by minors down the last column.  We obtain
$$\det(M-xI)=-x M_{n+9,n+9} + (-1)^nM_{1,n+9},\eqno(C.1)$$
where we use the notation $M_{i,j}$ for the $i,j$-minor of the matrix $M-xI$.  To evaluate $M_{n+9,n+9}$ and $M_{1,n+9}$, we expand again in minors along the last column to obtain
$$M_{n+9,n+9}=-x\det \hat m_1 + (-1)^n\det\hat m_2, \ \ \ \ \ M_{1,n+9}=\det\hat m_3,$$
where $\hat m_1=\pmatrix{A_1 & 0\cr
0&A_2(n)\cr}$, $\hat m_2 = \pmatrix{B_1& *\cr 0 & B_2(n)\cr}$, and $\hat m_3=\pmatrix{C_1 & *\cr 0 & C_2(n)\cr}$.  Here $A_1$, $B_1$, and $C_1$ do not depend on $n$, and $A_2(n)$, $B_2(n)$, and $C_2(n)$ are triangular matrices of size $(n-7)\times(n-7)$, $(n-8)\times(n-8)$ and $(n-8)\times(n-8)$ of the form
$$A_2(n)=\pmatrix{-x & 0&0\cr  *&\ddots&0\cr *& *& -x\cr}, \ \ B_2(n)=\pmatrix{1+x &* & *\cr   0&\ddots &* \cr 0 &0& 1+x\cr}, \ \ C_2(n)=\pmatrix{1&*&*\cr 0&\ddots&*\cr  0&0 & 1}.$$
Thus 
$$\det A_2(n)=(-x)^{n-7},\ \ \det B_2(n)=(1+x)^{n-8}, \ \ {\rm and\ }\det C_2(n)=1. \eqno(C.2)$$
 Since $A_1$, $B_1$, and $C_1$ do not depend on $n$, we may compute them using the matrix $M$ from the case $n=8$ to find 
 $$\det A_1=-x^6(x^8-x^5-x^3+1), \ \ \det B_1=-x^8,\ \  {\rm and} \ \det C_1=x^5+x^3-1.\eqno(C.3)$$
 Using (C.2) and (C.3) we find that the characteristic ploynomial of $M$ is equal to
 $$(-1)^n\left[ x^9(x+1)^{n-8}-x^{n+1}(x^8-x^5-x^3+1) +x^5+x^3-1\right]=(x^7-1)\chi_n(x)/(x^2-1),$$
 which completes the proof of $(i)$.  
 
 For the proof of $(ii)$, we use the ordered basis
 $$\{12,2a3,34, 45,  1g2, 2e3, 3c4, 04, 01, 03, 14, 25, 061, 072, 081, \dots, 0(2k-1)2,0(2k)1\},$$
 and for $(iii)$ we use the ordered basis 
 $$\{23, 3a4, 015, 13, 01, 12, 2d3, 3c4, 02, 04, 162, 073, 081, 192, \dots, 0(3k-1)1, 1(3k)2\}.  $$
 Otherwise, the proofs of cases $(ii)$ and $(iii)$ are similar.   We omit the details.  \qed

\bigskip\centerline{\bf References.}
\medskip

\item{[B]} A. Baker, The theory of linear forms in logarithms, in {\sl Transcendence Theory: Advances and Applications}, edited by A. Baker and D.W. Masser, Academic Press 1977, pp. 1--27.
\item{[BD1]}  E. Bedford and J. Diller,  Real and complex dynamics of a family of birational maps of the plane: the golden mean subshift. Amer. J. Math. 127 (2005), no. 3, 595--646.

\item{[BD2]} E. Bedford and J. Diller, Dynamics of a two parameter family of plane birational maps: maximal entropy, J. of Geometric Analysis, 16 (2006), no.\ 3, 409--430.   

\item{[BK]}  E. Bedford and KH Kim,   Periodicities in Linear Fractional Recurrences: Degree growth of birational surface maps, Mich. Math. J., 54 (2006), 647--670.  

\item{[BLS]} E. Bedford, M. Lyubich, and J. Smillie, Polynomial diffeomorphisms of $C\sp 2$. IV. The measure of maximal entropy and laminar currents. Invent. Math. 112 (1993), no. 1, 77--125.

\item{[BS1]} E. Bedford and J. Smillie, Polynomial diffeomorphisms of $C\sp 2$. II. Stable manifolds and recurrence. J. Amer. Math. Soc. 4 (1991), no. 4, 657--679. 

\item{[BS2]}  E. Bedford and J. Smillie,  Real polynomial diffeomorphisms with maximal entropy: Tangencies. Ann. of Math. (2) 160 (2004), no. 1, 1--26.

\item{[BC]} F. Beukers and R. Cushman, Zeeman's monotonicity conjecture, J. of Differential Equations, 143 (1998), 191--200.

\item{[C1]}  S. Cantat,  Dynamique des automorphismes des surfaces projectives complexes.  C. R. Acad. Sci. Paris S\'er. I Math. 328 (1999), no. 10, 901--906.

\item{[C2]}  S. Cantat,   Dynamique des automorphismes des surfaces $K3$.  Acta Math. 187 (2001), no. 1, 1--57.


\item{[DJS]}  J. Diller, D. Jackson, and A. Sommese,  Invariant curves for birational surface maps, Trans.\ AMS  359 (2007), 2973-2991.  math.AG/0505014

\item{[F]} W. Floyd,  Growth of planar Coxeter groups, P.V. numbers, and Salem numbers, Math.\ Ann.\ 293, 475--483 (1992).

\item{[FS]}  J.-E. Forn\ae ss and N. Sibony,  Classification of recurrent domains for some holomorphic maps. Math. Ann. 301 (1995), no. 4, 813--820. 

\item{[G]}   M. Gromov, Entropy, homology and semialgebraic geometry. S\'eminaire Bourbaki, Vol. 1985/86. Ast\'erisque No. 145-146 (1987), 5, 225--240. 

\item{[H]}   B. Harbourne, Rational surfaces with infinite automorpism group and no antipluricanonical curve, Proc. AMS, 99 (1987), 409--414.

\item{[KLR]} V.I. Kocic, G. Ladas, and I.W. Rodrigues, On rational recursive sequences, J. Math. Anal. Appl 173 (1993), 127-157.

\item{[KuL]}  M. Kulenovic and G. Ladas,  {\sl Dynamics of Second Order Rational Difference Equations}, CRC Press, 2002.


\item{[L]} R.C. Lyness, Notes 1581,1847, and 2952, Math. Gazette {\bf 26} (1942), 62, {\bf 29} (1945), 231, and {\bf 45} (1961), 201.

\item{[M1]}  C. McMullen, Dynamics on $K3$ surfaces: Salem numbers and Siegel disks. J. Reine Angew. Math. 545 (2002), 201--233. 

\item{[M2]}  C. McMullen,  Dynamics on blowups of the projective plane, Pub. Sci. IHES, 105 (2007), p. 49-89. 


\item{[P]} J. P\"oschel, On invariant manifolds of complex analytic mappings near fixed points. Exposition.  Math. 4 (1986), No. 2, 97--109.

\item{[U]} T. Ueda,  Critical orbits of holomorphic maps on projective spaces, J. of Geometric Analysis, 8 (1998), 319--334.

\item{[W]} L. Washington, {\sl Introduction to Cyclotomic Fields}, Springer-Verlag 1997.

\item{[Y]}  Y. Yomdin,  Volume growth and entropy. Israel J. Math. 57 (1987), no. 3, 285--300.

\item{[Z1]} C. Zeeman,  Geometric unfolding of a difference equation, lecture.

\item{[Z2]}  E. Zehnder, A simple proof of a generalization of a theorem by C. L. Siegel. Geometry and topology (Proc. III Latin Amer. School of Math., Inst. Mat. Pura Aplicada CNPq, Rio de Janeiro, 1976), pp. 855--866. Lecture Notes in Math., Vol. 597, Springer, Berlin, 1977. 

\item{[Zh]} D-Q Zhang,  Automorphism groups and anti-canonical curves, Math.\ Res.\ Lett.\ 15 (2008), 163--183.

\bigskip
\rightline{Indiana University}

\rightline{Bloomington, IN 47405}

\rightline{\tt bedford@indiana.edu}

\bigskip
\rightline{Florida State University}

\rightline{Tallahassee, FL 32306}

\rightline{\tt kim@math.fsu.edu}

\bye